\documentclass[a4paper,11pt]{article}
\usepackage{latexsym,amsfonts,amsmath}
\usepackage{amssymb}
\usepackage{color,soul}
\usepackage{xcolor}
\usepackage[mathscr]{eucal}
\usepackage{enumerate}
\usepackage{enumitem}

\usepackage{pgf}
\usepackage{tikz}
\usetikzlibrary{arrows,automata}

\usepackage{fullpage}

\definecolor{Dgreen}{rgb}{0,0.5,0}

\def\tto{\twoheadrightarrow}

\def\1{\mathbf{1}}
\def\0{\mathbf{0}}

\def\NN{\mathbb{N}}

\def\RR{\mathbb{R}}

\def\car{{\mathcal{C}ar}}

\def\XX{\mathbf{X}}
\def\AA{\mathbf{A}}
\def\KK{\mathbf{K}}
\def\SS{\mathbf{S}}

\newcommand{\bcal}[1]{\boldsymbol{\mathcal{#1}}}

\newcommand{\bfrak}[1]{\boldsymbol{\mathfrak{#1}}}

\def\argmax{\mathop{\hbox{\rm arg max}}\limits}

\def\limsup{\mathop{\overline{\lim}}}
\def\ds{\displaystyle}

\def\wstar{\ {\stackrel{*}{\rightharpoonup}\ }}

\definecolor{Dgreen}{rgb}{0,0.5,0}

\newcounter{hypot}

    \newenvironment{hypot}{\begin{list}
      {\hspace{\labelsep}\bfseries Assumption \Alph{hypot}.}
      {\leftmargin=0pt
       \labelwidth=0cm
       \refstepcounter{hypot}
       \def\makelabel##1{##1}}}{\end{list}}

\newcounter{assump}

\newenvironment{assump}{\begin{list}
      {\hspace{\labelsep} (\Alph{hypot}.\arabic{assump})}
      {\leftmargin=0pt
       \labelwidth=0cm
       \usecounter{assump}
       }}{\end{list}}

\newcounter{cond}

\begin{document}
\newtheorem{theorem}{Theorem}[section]
\newtheorem{proposition}[theorem]{Proposition}
\newtheorem{lemma}[theorem]{Lemma}
\newtheorem{corollary}[theorem]{Corollary}
\newtheorem{definition}[theorem]{Definition}
\newtheorem{remark}[theorem]{Remark}
\newtheorem{conjecture}[theorem]{Conjecture}
\newtheorem{assumption}[theorem]{Assumption}

\bibliographystyle{plain}

\title{Stationary Markov Nash equilibria for nonzero-sum constrained ARAT Markov games}

\author{Fran\c{c}ois Dufour\footnote{
Institut Polytechnique de Bordeaux; INRIA Bordeaux Sud Ouest, Team: ASTRAL; IMB, Institut de Math\'ematiques de Bordeaux, Universit\'e de Bordeaux, France \tt{francois.dufour@math.u-bordeaux.fr}}
\and
Tom\'as Prieto-Rumeau\footnote{Statistics Department, UNED, Madrid, Spain. e-mail: {\tt{tprieto@ccia.uned.es}}\quad {(Author for correspondence)}}}
\date{}
\maketitle
\begin{abstract}
We consider a nonzero-sum Markov game on an abstract measurable state space with compact metric action spaces. The goal of each player is to maximize his respective discounted payoff function under the condition that some constraints on a discounted payoff are satisfied. We are interested in the existence of a Nash or noncooperative equilibrium. Under suitable conditions, which include absolute continuity of the transitions with respect to some reference probability measure, additivity of the payoffs and the transition probabilities (ARAT condition), and  continuity  in action of the payoff functions and the density function of the transitions of the system, we establish the existence of a constrained stationary Markov Nash equilibrium,
that is, the existence of stationary Markov strategies for each of the players yielding an optimal profile within the class of all history-dependent profiles.
\end{abstract}
{\small 
}

\section{Introduction}

Nonzero-sum stochastic games are nowadays a largely developed and yet very active field of research. It started in the 1950's with the pioneering work of J. Nash and it has developed in many directions. It is impossible to give a complete overview here, though
we refer the reader to the recent comprehensive survey \cite{jaskiewicz18} and the citations therein.
The existence of stationary equilibria is a central problem in the theory nonzero-sum stochastic games, and it can be traced back to \cite{himmelberg76}. 

Results on the existence
of stationary equilibria have then been obtained for different families of games.
We can mention  \cite{himmelberg76,nowak87,parthasarathy82} for ARAT (Additive Reward and Additive Transition) games, \cite{parthasarathy89} for games  with finite action spaces and state independent transition kernel, \cite{nowak92} for  stationary correlated equilibria (see  \cite{duffie94} for related results), \cite{duggan12} for noisy stochastic games, \cite{nowak03} for a class of games where the transition probability measure is a convex combination of a finite number of probability measures and where the coefficients depend on the state and actions variables, and \cite{he-sun17} for a family of games satisfying a general condition called ``decomposable coarser transition kernels''.

For the particular case  of discounted stochastic ARAT games, its analysis goes back to the 1970s and 1980s with results on the existence of stationary equilibria under the assumption that the state space is Borel and the action spaces are finite \cite{himmelberg76,parthasarathy82}. These first results were then generalized in \cite{nowak87} in several directions by considering an ARAT game with general measurable state space and compact action spaces, showing the existence of nonrandomized stationary
$\epsilon$-equilibria assuming that the initial distribution is nonatomic. In \cite{jasmiewicz15}, the authors assume that the state space is Borel, the action space is finite, and also that the transition probability distribution is absolutely continuous with respect to a nonatomic probability measure, and then show the existence of pure Nash equilibria in the class of stationary almost Markov strategies. Other type of criteria are analyzed in \cite{kuenle99,thuijsman97} for ARAT stochastic games. More specifically, the expected total cost criterion is studied in \cite{kuenle99}, while the authors focus on the limiting average criterion in \cite{thuijsman97}.

In the aforementioned references on \textit{unconstrained} games, the usual technique is to reduce the infinite horizon discounted game to a one-shot game with some terminal payoff function $v$.
Each player addresses his associated control problem (with terminal payoff $v$) by solving the corresponding \textit{dynamic programming equation}. 
The main idea is to obtain a terminal payoff function~$v^*$ satisfying a suitable fixed point property (related to fixed point theorems  for correspondences). Nash equilibria for the players are then obtained from this terminal payoff function $v^*$. To use such fixed point results, a \textit{convexification} step is needed and a final \textit{de-convexification} step becomes also necessary in order to get back to the original problem.
This important procedure, obtained in \cite{nowak92}, is a key point to prove the existence of stationary equilibria for different types of games  or different notions of equilibrium
 (see \cite{nowak07} for examples arising from economic theory or  \cite{nowak92} for  stationary equilibria with public randomization). 
Finally, \cite{he-sun17} gives a general de-convexification procedure based on the already mentioned notion  of a decomposable coarser transition kernel.

\bigskip

A natural extension of the \textit{unconstrained} game models described above is to impose constraints on the players. More precisely, each player has a constraint function and he has to maximize his own discounted payoff under the condition that the expected discounted payoff associated to the constraint function is above some given level. A Nash equilibrium consists of policies of the players such that the constraints of each player are satisfied and, in addition, no player improves his payoff when unilaterally varying his strategy while still satisfying his own constraint.  Constrained games are, by far, less developed than the unconstrained counterpart.  

In the reference \cite{Altman-Shwartz-96}, the authors study a constrained game with finite state and actions spaces. 
In the same vein as for constrained Markov decision processes, the approach developed in \cite{Altman-Shwartz-96} consists in considering the family of occupation measures of the state-action process of each player. From these occupation measures, the corresponding policies of the players are identified by disintegration. Then, when fixing the policies of the remaining players, each player finds his optimal constrained response by solving a linear programming problem stated in, again, the family of occupation measures. In this way, a correspondence on the family of occupation measures is defined, and any fixed point of this correspondence yields a constrained Nash equilibrium.
In \cite{Mena-OHL-constrained-games}, these results are generalized to a game model with countable state space and compact action spaces which is, roughly speaking, ``nearly finite'' in each transition of the dynamic system \cite[Assumption~3.4]{Mena-OHL-constrained-games}.  Additional results have been obtained in, e.g.,  \cite{Altman-Avrachenkov-08,Singh14}, where the authors consider a finite state and action game model with independent state and action dynamics, in which interaction between players occurs only on the payoff functions. 

In this paper,  we will consider a constrained nonzero-sum stochastic game with ARAT structure under the discounted payoff criterion. The state space is an abstract measurable space, the action spaces of the players are compact metric spaces, the payoff and constraint functions are bounded Carath\'eodory (i.e., measurable in state, continuous in action) functions. It is important to emphasize that the results presented in \cite{Altman-Shwartz-96} cannot be generalized to the general state space case.
Indeed, compared to \cite{Altman-Shwartz-96}, dealing with constraints and extending the state space from a finite or countable space to a general measurable (not even metric) space entails serious technical difficulties and it is ---we believe--- far from being straightforward.
To work in such a general context, one needs to introduce additional hypotheses including continuity type assumptions of the transition kernel and  the ARAT separability condition similar to those used in \cite{nowak87} for unconstrained games.
In particular, the ARAT condition is used to uniquely identify a Markov policy of a player starting from a state-action occupation measure and to characterize the corresponding payoffs.  We use the weak-strong topology for the space of occupation measures, which is compatible with the so-called narrow topology (also known as the stable topology) of Young measures, and the latter will be used for the spaces of Markov policies of the players.
In this fairly general setting, we will be able to establish the existence of a constrained Nash equilibrium ---consisting of stationary Markov policies of the players--- within the class of all history-dependent policies of the players.
Finally, it is worth mentioning that the
existence of a stationary Markov equilibrium is an important and difficult issue in the literature on stochastic games. In support of this, we can cite A. Nowak who wrote in \cite{nowak85}: \textit{This question is highly nontrivial and not as yet completely resolved}. As emphasized in \cite{he-sun17}, this issue still remains an important problem.
One aspect of our contribution is to show that this type of existence result holds for a class of constrained stochastic games with general state space.
 
 \bigskip
 
The rest of the paper is organized as follows. In Section \ref{section2} we define the constrained game model and state our assumptions. 
We introduce the occupation measures of the players in Section~\ref{section3} and we also state some preliminary results. Finally, the existence of constrained Nash equilibria is addressed in Section \ref{section4}. We note that we deal here 
with the two-player game so as to handle a simpler notation. Our results  can be easily generalized to an $m$-player model.

\section{Model and assumptions}\label{section2}
\setcounter{equation}{0}

\subsection{Notation and preliminary results}

On a measurable space $(\mathbf{\Omega},\mathcal{F})$ we will consider the set of finite signed measures $\bcal{M}(\mathbf{\Omega})$, the set of finite nonnegative measures $\bcal{M}^+(\mathbf{\Omega})$, and the set of probability measures~$\bcal{P}(\mathbf{\Omega})$.
On $\bcal{P}(\mathbf{\Omega})$, the $s$-topology is the coarsest topology that makes 
$\mu\mapsto \mu(D)$
continuous for every $D\in\mathcal{F}$.
For $\lambda\in\bcal{P}(\mathbf{\Omega})$, let us consider the set $\bcal{P}_{\lambda}(\mathbf{\Omega})$ of probability measures~$\eta$ on~$\mathbf{\Omega}$ which are absolutely continuous with respect to $\lambda$; in symbols, $\eta\ll\lambda$.
The set $\bcal{P}^{e}_{\lambda}(\mathbf{\Omega})$ consists of all probability measures $\eta$ on $\mathbf{\Omega}$  which are equivalent to $\lambda$, that is: $\eta\ll\lambda$ and $\lambda\ll\eta$, written $\eta\sim\lambda$.

For a product of measurable spaces we will always consider the product $\sigma$-algebra. 
Let $(\mathbf{\Omega},\mathcal{F})$ and $(\mathbf{\Omega'},\mathcal{F'})$ be two measurable spaces. A kernel on $\mathbf{\Omega'}$ given $\mathbf{\Omega}$ is a mapping $Q:\mathbf{\Omega}\times\mathcal{F'}\rightarrow\RR^+$ such that $\omega\mapsto Q(B|\omega)$ is measurable on $(\mathbf{\Omega},\mathcal{F})$ for every $B\in\mathcal{F}'$,   and  $B\mapsto Q(B|\omega)$ is in $\bcal{M}^+(\mathbf{\Omega'})$ for every $\omega\in\mathbf{\Omega}$. If $Q(\mathbf{\Omega'}|\omega)=1$ for all $\omega\in\mathbf{\Omega}$ then we say that $Q$ is a \textit{stochastic kernel} and if $Q(\mathbf{\Omega'}|\omega)\le1$ for all $\omega\in\mathbf{\Omega}$ then we say that $Q$ is a \textit{substochastic kernel}.
Let $Q$ be a kernel on $\mathbf{\Omega'}$ given $\mathbf{\Omega}$ and let $f$ be a bounded measurable function $f:\mathbf{\Omega'}\rightarrow\RR$. We will denote by $Qf:\mathbf{\Omega}\rightarrow\RR$ the measurable function
$$Qf(\omega)=\int_\mathbf{\Omega'} f(z)Q(dz|\omega)\quad\hbox{for $\omega\in\mathbf{\Omega}$}.$$
If $Q$ is a stochastic 
(or substochastic) 
kernel on $\mathbf{\Omega'}$ given $\mathbf{\Omega}$ and  $\mu\in\bcal{P}(\mathbf{\Omega})$, we denote by $\mu Q$ the probability measure (or finite measure) 
 on $(\mathbf{\Omega'},\mathcal{F}')$ given by
$$B\mapsto \mu Q\,(B)= \int_{\mathbf{\Omega}} Q(B|\omega) \mu(d\omega)\quad\hbox{for $B\in\mathcal{F}'$}.$$
We will also write $\mu Q\,(d\omega')=\int_{\mathbf{\Omega}}Q(d\omega'|\omega)\mu(d\omega)$.
In addition, we define $\mu\otimes Q$ as the probability measure (or finite measure) on the product space $(\mathbf{\Omega}\times\mathbf{\Omega'},\mathcal{F}\otimes\mathcal{F}')$ given by 
$$(\mu\otimes Q)( A\times B)= \int_A Q(B|\omega)\mu(d\omega)\quad\hbox{for $(A,B)\in\mathcal{F}\times\mathcal{F}'$}
.$$
We will write $(\mu\otimes Q)(d\omega,d\omega')=Q(d\omega'|\omega)\mu(d\omega)$. 
Given $\mu\in\bcal{M}(\mathbf{\Omega}\times\mathbf{\Omega}')$,  the  marginal measures are $\mu^{\mathbf{\Omega}}\in\bcal{M}(\mathbf{\Omega})$ and $\mu^{\mathbf{\Omega}'}\in\bcal{M}(\mathbf{\Omega'})$ defined by 
$\mu^{\mathbf{\Omega}}(\cdot)=\mu(\cdot\times \mathbf{\Omega}')$ and
$\mu^{\mathbf{\Omega}'}(\cdot)=\mu(\mathbf{\Omega}\times\cdot)$.

Throughout this paper, any metric space $\SS$ will be endowed with its Borel $\sigma$-algebra $\bfrak{B}(\SS)$. Also, when considering the product of a finite family of metric spaces, we will consider the product topology (which makes the product again a metric space).
We say that $f:\mathbf{\Omega}\times\SS\rightarrow\SS'$, where~$\SS'$ is a metric space, is a \textit{Carath\'eodory function} if $f(\cdot,s)$ is measurable on~$\mathbf{\Omega}$ for every $s\in\SS$ and $f(\omega,\cdot)$ is continuous on $\SS$ for every $\omega\in \mathbf{\Omega}$. The family of the so-defined Carath\'eodory functions is denoted by $\car(\mathbf{\Omega}\times\SS,\SS')$. The family of Carath\'eodory functions which, in addition, are bounded is denoted by $\car_b(\mathbf{\Omega}\times\SS,\SS')$. 
When the metric space~$\SS$ is separable then any $f\in\car(\mathbf{\Omega}\times\SS,\SS')$ is a jointly measurable function on $(\mathbf{\Omega}\times\mathbf{S},\mathcal{F}\otimes\bfrak{B}(\mathbf{S}))$; see \cite[Lemma 4.51]{aliprantis06}.

If $\SS$ is a Polish space (a complete and separable metric space), on $\bcal{M}(\mathbf{\Omega}\times\SS)$ we will consider the  $ws$-topology (weak-strong topology) which is the coarsest topology for which the mappings
$$
\mu\mapsto \int_{\mathbf{\Omega}\times\SS} f(\omega,s)\mu(d\omega,ds)
$$
for $f\in\car_b(\mathbf{\Omega}\times\SS,\RR)$ are continuous.
There are other equivalent definitions for this topology as discussed, for instance, in \cite{balder01,jacod81,schal75}.
Note that $\bcal{P}(\mathbf{\Omega}\times\SS)$ is a closed subset of $\bcal{M}(\mathbf{\Omega}\times\SS)$.
It is clear that convergence in the $ws$-topology in $\bcal{P}(\mathbf{\Omega}\times\SS)$  implies convergence of the $\mathbf{\Omega}$-marginal probability measures in the $s$-topology.

Inequality $\ge$ in $\RR^p$ means a componentwise inequality $\ge$, while the inequality $>$ in $\RR^p$ is a componentwise strict inequality $>$. Let $\1_{p}\in\RR^p$ be the vector with all components equal to one.

The next result characterizes relative compactness in  $\bcal{P}(\mathbf{\Omega}\times\SS)$. It ressembles to \cite[Theorem~3.10]{schal75} except that, here, $\mathbf{\Omega}$ is a measurable space, whereas  \cite{schal75}  assumes that $\mathbf{\Omega}$ is Borel.

\begin{proposition}
\label{Criterion-relatively-ws-compactness}
Let $(\mathbf{\Omega},\mathcal{F})$ be a measurable space and let $\SS$ be a Polish space, and consider the
$ws$-topology on $\bcal{P}(\mathbf{\Omega}\times\SS)$. A necessary and sufficient condition for a set $K\subseteq \bcal{P}(\mathbf{\Omega}\times\SS)$ to be relatively compact is that   for any sequence $\{f_n\}_{n\in\NN}$ of functions  $f_n\in\car_b(\mathbf{\Omega}\times\SS,\RR)$ which decrease pointwise to zero (i.e., $f_n(\omega,s)\downarrow0$ for each $(\omega,s)\in\mathbf{\Omega}\times\SS$) we have
\begin{equation}
\lim_{n\rightarrow\infty} \sup_{\mu\in K} \int_{\mathbf{\Omega}\times\SS} f_{n} d\mu = 0.
\label{uniform-principle-relatively-compact}
\end{equation}
\end{proposition}
\textbf{Proof.} The $ws$-topology being the weak topology induced by bounded Carath\'eodory functions, we use
 \cite[Lemma 3.1]{schal75} for the necessary condition. For sufficiency, by  \cite[Theorem 5.2]{balder01} we have that $K$ is relatively compact if and only if 
$K^{\mathbf{\Omega}}=\{\mu^{\mathbf{\Omega}} : \mu\in K\}$ and $ K^\SS=\{\mu^{\SS} : \mu\in K\}$
are relatively $s$-compact (for the $s$-topology) and relatively $w$-compact (for  weak convergence of probability measures on $\SS$), respectively. For relative $s$-compactness, we use \cite[Lemma 3.5]{schal75}.
If $\{g_n\}$ is a sequence of bounded measurable functions on $\mathbf{\Omega}$ decreasing pointwise to zero, we have
$\sup_{\mu\in K} \int_{\mathbf{\Omega}} g_n d\mu^{\mathbf{\Omega}}=\sup_{\mu\in K}\int_{\mathbf{\Omega}\times\SS} \bar{g}_nd\mu\rightarrow0$
by hypothesis, where $ \bar{g}_n(\omega,s)=g_n(\omega)$,
showing that  $K^{\mathbf{\Omega}}$ is relatively $s$-compact. 
Similarly, if $\{h_n\}$ is a sequence of bounded continuous functions on $\SS$ which decrease pointwise to zero then we have
$\sup_{\mu\in K} \int_{\SS} h_n d\mu^{\SS}=\sup_{\mu\in K}\int_{\mathbf{\Omega}\times\SS} \bar{h}_nd\mu\rightarrow0$
where $ \bar{h}_n(\omega,s)=h_n(s)$.
By \cite[Lemma 3.2]{schal75}, we obtain that  $K^\SS$ is  relatively $w$-compact.
\hfill$\Box$
\\[10pt]\indent
We now state  without proof  two results that we will use repeatedly in the following.
\begin{proposition}\label{prop-locally-convex}
If $(\mathbf{\Omega},\mathcal{F})$ is a  measurable space and $\SS$ is a Polish space, then $\bcal{M}(\mathbf{\Omega}\times\SS)$ equipped with the $ws$-topology
is a locally convex Hausdorff topological vector space.
\end{proposition}

\begin{lemma}[Disintegration lemma]
\label{lemma-disintegration} 
Let $(\mathbf{\Omega},\mathcal{F})$ be a measurable space and let $\SS$ be a Polish space. Let $\varphi:\mathbf{\Omega}\tto\SS$ be a weakly measurable correspondence with nonempty closed values, and let~$\mathbf{K}$ be the graph of the correspondence. For every  $\mu\in\bcal{P}(\mathbf{\Omega}\times\SS)$ such that $\mu(\mathbf{K})=1$ there exists a stochastic kernel $Q$ on $\SS$ given $\mathbf{\Omega}$ such that 
\begin{equation}\label{eq-dudley-product}
\mu= \mu^{\mathbf{\Omega}}\otimes Q 
\end{equation}
and such that
$Q(\varphi(\omega)|\omega)=1$ for each $\omega\in\mathbf{\Omega}$. Moreover, $Q$ is unique $\mu^{\mathbf{\Omega}}$-almost surely, meaning that if $Q$ and $Q'$ are two stochastic kernels that satisfy \eqref{eq-dudley-product} then for all $\omega$ in a set of $\mu^{\mathbf{\Omega}}$-probability one, the probability measures $Q(\cdot|\omega)$ and $Q'(\cdot|\omega)$ coincide.
\end{lemma}

\subsection{Game model and assumptions}

\paragraph{Elements of the game model.}
The elements of the game model are the following.
 \begin{enumerate}
\item [(a)] $(\XX,\bfrak{X})$ is an abstract measurable space, where $\XX$ stands for the state space of the system.
\item [(b)] The separable metric spaces $\mathbf{A}_{1}$ and $\AA_2$ represent the action sets for player 1 and 2, respectively. For each $x\in\XX$ the nonempty measurable sets $\AA_1(x)\subseteq\AA_1$ and $\AA_2(x)\subseteq\AA_2$ 
are the set of feasible actions for players 1 and 2, respectively, when the system is in state $x\in\XX$.
Let   $\AA=\AA_1\times\AA_2$ and 
$\mathbf{A}(x)=\mathbf{A}_{1}(x)\times \mathbf{A}_{2}(x)$ for each $x\in \mathbf{X}$.
\item [(c)] Fix a player $i=1,2$. The measurable  functions $r_{i}:\XX\times\AA\rightarrow\RR$ and $c_i:\XX\times\AA\rightarrow\RR^p$ are the reward  and the constraint functions, respectively, of  player $i$. The constraint constant for player~$i$ is $\rho_i\in\RR^p$. Here, $p\ge1$ is a fixed integer.
Write $\rho=(\rho_1,\rho_2)\in\RR^{p}\times\RR^p$.
\item [(e)] The transitions of the system are given by a stochastic kernel $Q$ on $\mathbf{X}$ given $\mathbf{X}\times\mathbf{A}$. 
\item [(f)] The discount factor is $\beta\in (0,1)$.
\item [(g)] The initial distribution is the probability measure $\eta\in\bcal{P}(\XX)$.
\end{enumerate}
This game model will be denoted by $\mathcal{G}(\eta,\rho)$. The reason for this notation is that, 
in the forthcoming, we will need  the initial distribution~$\eta$ and the constraint constants $(\rho_1,\rho_2)$  to vary, while the other elements of the game model  will remain fixed.
\paragraph{Policies of the players.}
We define the sets 
$\mathbf{H}_0=\XX$ and $
\mathbf{H}_t=(\XX\times\AA_1\times\AA_2)^t\times\XX$
for $t\ge1$. The elements of $\mathbf{H}_t$ will be usually denoted by $(x_0,a_0,b_0,\ldots,x_{t-1},a_{t-1},b_{t-1},x_t)$. 
\begin{definition}\label{def-policies}
\begin{itemize}
\item[(i).] A  policy for player $i$ (with $i=1,2$) is a sequence $\{\pi_i^t\}_{t\in\NN}$ of stochastic kernels on $\AA_i$ given $\mathbf{H}_t$ such that 
$$\pi_i^t(\AA_i(x_t)|x_0,a_0,b_0,\ldots,x_t)=1\quad\hbox{for all $t\ge0$ and $h_t=(x_0,a_0,b_0,\ldots,x_t)\in\mathbf{H}_t$}.$$
The family of all policies of player $i$ is denoted by $\mathbf{\Pi}_i$. 
\item[(ii).] Let $\mathbf{M}_i$ be the family of stochastic kernels $\pi_i$ on $\AA_i$ given $\XX$ such that $\pi_i(\AA_i(x)|x)=1$ for each $x\in\XX$.
We say that $\{\pi_i^t\}_{t\in\NN}\in\mathbf{\Pi}_i$ is a stationary Markov policy for player $i$ (with $i=1,2$) if there is some $\pi_i\in\mathbf{M}_i$ which verifies
$$\pi_i^t(\cdot|x_0,a_0,b_0,\ldots,x_t)=\pi_i(\cdot|x_t)\quad\hbox{for all $t\ge0$ and $h_t=(x_0,a_0,b_0,\ldots,x_t)\in\mathbf{H}_t$}.$$
\end{itemize}
\end{definition}

The class $\mathbf{\Pi}_i$ is also referred to as the family of
\textit{history-dependent policies} for player $i$. 
We can identify the class of stationary Markov  policies for player $i$ with  $\mathbf{M}_i$, and so $\mathbf{M}_i\subseteq\mathbf{\Pi}_i$.
In the sequel, we will refer to $\mathbf{M}_i$ as to the family of \textit{stationary Markov policies}  for player $i$. 
Our conditions below will ensure that these classes of policies are nonempty. 

\paragraph{Construction of the state-actions process.}
We address the existence of a probability space supporting the dynamic system of the game model $\mathcal{G}(\eta,\rho)$.
On the product canonical space $\mathbf{H}_\infty=(\XX\times\AA_1\times\AA_2)^\NN$, endowed with its product $\sigma$-algebra $(\bfrak{X}\otimes\bfrak{B}(\AA_1)\otimes\bfrak{B}(\AA_2))^\NN$, define $(X_t,A_t,B_t)_{t\in\NN}$ as the coordinates projection operators and write  $H_t=(X_0,A_0,B_0,\ldots,X_t)$.
Given an initial distribution $\eta\in\bcal{P}(\XX)$  and policies $(\pi_1,\pi_2)\in\mathbf{\Pi}_1\times\mathbf{\Pi}_2$,  we can uniquely define a probability measure
$\mathbb{P}_{\eta,\pi_1,\pi_2}$ on $\mathbf{H}_\infty$ which verifies the following properties for any $D\in\bfrak{X}$, $D_1\in\bfrak{B}(\AA_1)$, and $D_2\in\bfrak{B}(\AA_2)$, and every $t\ge0$.  Firstly,
$
\mathbb{P}_{\eta,\pi_1,\pi_2}\{X_0\in D\}=\eta(D)$.
Secondly,
$$
\mathbb{P}_{\eta,\pi_1,\pi_2}(A_t\in D_1,B_t\in D_2 |H_t)=\pi^t_1(D_1|H_t)\cdot\pi^t_2(D_2|H_t).
$$
is  the conditional distribution of the actions.
Finally, distribution of  the next state of the system is 
$\mathbb{P}_{\eta,\pi_1,\pi_2}(X_{t+1}\in D| X_0,A_0,B_0,\ldots,X_t,A_t,B_t)=Q(D|X_t,A_t,B_t)$.
The expectation operator associated to $\mathbb{P}_{\eta,\pi_1,\pi_2}$ will be denoted by $\mathbb{E}_{\eta,\pi_1,\pi_2}$.

\paragraph{Payoffs of the players and equilibrium.}
Our assumptions  below will ensure that the following expressions are well defined and finite. 
Consider the game model $\mathcal{G}(\eta,\rho)$ and suppose the players use the  policies $(\pi_1,\pi_2)\in\mathbf{\Pi}_1\times\mathbf{\Pi}_2$. The total expected discounted reward for player $i$ is
$$R_i(\eta,\pi_1,\pi_2)=(1-\beta)\mathbb{E}_{\eta,\pi_1,\pi_2} \Big[ \sum_{t=0}^\infty \beta^t r_i(X_t,A_t,B_t)\Big]\in\RR,$$
while the total expected discounted constraint for player $i$ is 
$$C_i(\eta,\pi_1,\pi_2)=(1-\beta)\mathbb{E}_{\eta,\pi_1,\pi_2} \Big[ \sum_{t=0}^\infty \beta^t c_i(X_t,A_t,B_t)\Big]\in\RR^p.$$
We say that the policies $(\pi_1,\pi_2)\in\mathbf{\Pi}_1\times\mathbf{\Pi}_2$ satisfy the constraint of player $i$  when $C_i(\eta,\pi_1,\pi_2)\ge\rho_i$.

\begin{definition}\label{def-nash} Consider the game model $\mathcal{G}(\eta,\rho)$.
The pair of policies $(\pi^*_1,\pi^*_2)\in\mathbf{\Pi}_1\times\mathbf{\Pi}_2$ is a  constrained Nash equilibrium in the class of all history-dependent policies
$\mathbf{\Pi}_1\times\mathbf{\Pi}_2$ when:
\begin{itemize}
\item[(i)] The constraints of both players are satisfied:
$C_i(\eta,\pi^*_1,\pi^*_2)\ge\rho_i$ for $i=1,2$.
\item[(ii)]   The following conditions hold:
$$\forall\pi_1\in\mathbf{\Pi}_1,\ C_1(\eta,\pi_1,\pi^*_2)\ge\rho_1 \ \Rightarrow \ R_1(\eta,\pi^*_1,\pi^*_2)\ge R_1(\eta,\pi_1,\pi^*_2).$$
$$\forall\pi_2\in\mathbf{\Pi}_2,\ C_2(\eta,\pi^*_1,\pi_2)\ge\rho_2 \ \Rightarrow \ R_2(\eta,\pi^*_1,\pi^*_2)\ge R_2(\eta,\pi^*_1,\pi_2).$$
\end{itemize}
\end{definition}

\paragraph{Assumptions on the game model.}
First of all, we introduce a Slater-type condition. 

\begin{definition}\label{def-slater}
We say that the game model $\mathcal{G}(\eta,\rho)$
satisfies the Slater condition if
$$\forall\,\pi_1\in\mathbf{M}_1, \exists\, \pi'_2\in\mathbf{M}_2, C_2(\eta,\pi_1,\pi'_2)>\rho_2\quad\hbox{and}\quad
\forall\,\pi_2\in\mathbf{M}_2, \exists\, \pi'_1\in\mathbf{M}_1, C_1(\eta,\pi'_1,\pi_2)>\rho_1.$$
\end{definition}

An equivalent formulation of this condition is discussed in
 Remark \ref{rem-Markov-sufficient}.
This is the usual form of the Slater condition for constrained games; see \cite[Assumption $\mathit{\Pi_2}$]{Altman-Shwartz-96}  and \cite[Assumption 3.3.c]{Mena-OHL-constrained-games}. A stronger uniform Slater condition is imposed in \cite[Assumption $\mathit{\Pi_2}$]{Altman-Avrachenkov-08} and in \cite[Assumption A.2]{Singh14}.

Next we state our conditions on the game model. Notice that Assumptions \ref{Assumption-reward}--\ref{Assumption-transition-Q} impose precisely the additive reward additive transition (ARAT) character of the game. 
\begin{hypot}
\label{Assump}
\item \mbox{ } For an initial distribution $\nu\in\bcal{P}(\XX)$ and constraint constants $\theta=(\theta_1,\theta_2)$ in $\RR^p\times\RR^p$, we suppose that the game model $\mathcal{G}(\nu,\theta)$ satisfies the following conditions.
\begin{assump}
\item \label{Assumption-state} The $\sigma$-algebra $\bfrak{X}$ is countably generated.
\item \label{Assumption-control} For each $i=1,2$ the action set $\mathbf{A}_{i}$ is compact 
and the correspondence from $\mathbf{X}$ to $\mathbf{A}_{i}$ defined by $x\mapsto \mathbf{A}_{i}(x)$ is weakly measurable with nonempty compact values.
 \item \label{Assumption-reward} For each player $i$ (with $i=1,2$)  there exist Carath\'eodory functions  $r_i^j\in\car_b(\XX\times\AA_j,\RR)$ and  $c_i^j\in\car_b(\XX\times\AA_{j},\RR^p)$ for $j=1,2$ such that, for every $(x,a_1,a_2)\in\XX\times \AA$,
 $$r_{i}(x,a_1,a_2)=r_{i}^{1}(x,a_1)+r_{i}^{2}(x,a_2) \quad\hbox{and}\quad
c_{i}(x,a_1,a_2)=c_{i}^{1}(x,a_1)+c_{i}^{2}(x,a_2).$$
\item \label{Assumption-transition-Q} There exist a probability measure $\lambda\in\bcal{P}(\XX)$ and measurable functions $q_i: \XX\times\XX\times\AA_{i}\rightarrow\RR^+$ for $i=1,2$ such that  for every $D\in\bfrak{X}$  and every $(x,a_1,a_2)\in\XX\times\AA$
$$Q(D|x,a_1,a_2)=\int_D [q_{1}(y,x,a_{1})+ q_{2}(y,x,a_{2}) ]\lambda(dy).$$
Also, the following continuity condition holds: 
for each $i=1,2$ and every $x\in\XX$ we have
$$\lim_{n\rightarrow\infty} \int_{\XX} |q_{i}(y,x,b_{n})-q_{i}(y,x,b)| \lambda(dy)=0$$
whenever $b_n\rightarrow b$ in $\AA_i$. 
Moreover, the probability space
 $(\XX,\bfrak{X},\lambda)$ is complete. 
 \item
\label{Assump-Slater} The initial distribution $\nu\in\bcal{P}(\XX)$ and the constraint constants $(\theta_{1},\theta_{2})\in\RR^{p}\times\RR^{p}$ satisfy the Slater condition in Definition \ref{def-slater}.
\end{assump}
\end{hypot}

\textit{
In the sequel we will suppose that  Assumption \ref{Assump} holds  with no need of an explicit mention. }
Our main result in this paper ---Theorem \ref{th-main} below--- establishes that,  for the game model $\mathcal{G}(\nu,\theta)$, there exists a pair of stationary Markov policies $(\pi^*_1,\pi^*_2)\in\mathbf{M}_1\times\mathbf{M}_2$ which is a constrained Nash equilibrium within the class $\mathbf{\Pi}_1\times\mathbf{\Pi}_2$ of all history-dependent policies.
Some comments on Assumption \ref{Assump} are given next.
\begin{remark}\label{rem-after-assumption}
\begin{enumerate}
\item[(a).]  There is no loss of generality in assuming that the initial distribution $\nu$  is absolutely continuous with respect to $\lambda$. Indeed, we can replace $\lambda$ in Assumption \ref{Assump} with  $\bar{\lambda}=(\nu+\lambda)/2$, so that $\nu\in\bcal{P}_{\bar\lambda}(\XX)$.
It is then easily seen that the functions $\bar{q}_{i}(y,x,a_{i})=q_{i}(y,x,a_{i})\frac{d\lambda}{d\bar{\lambda}}(y)$  satisfy Assumption \ref{Assumption-transition-Q} and that $(\XX,\bfrak{X},\bar{\lambda})$ is  a complete probability space as well. 
Hence, from now on we will consider that $\nu\in\bcal{P}_\lambda(\XX)$.

\item[(b).] By Assumption \ref{Assumption-control}, the correspondences $x\mapsto \AA_1(x)$ and $x\mapsto\AA_2(x)$ are measurable \cite[Lemma 18.2]{aliprantis06} and they have measurable graph  \cite[Theorem 18.6]{aliprantis06}. Therefore,  
\begin{eqnarray*}
 \KK_{i}&=&\{(x,a_{i}):\XX\times\AA_{i} : a_{i}\in\AA_{i}(x) \}\in \bfrak{X}\otimes\bfrak{B}(\AA_{i}) \quad\hbox{for $i=1,2$}\\
 \KK\hphantom{_i}&=&\{(x,a_{1},a_{2}):\XX\times\AA_{1}\times\AA_{2} : (a_1,a_2)\in\AA(x)\}\in \bfrak{X}\otimes\bfrak{B}(\AA_1)\otimes\bfrak{B}(\AA_2).
\end{eqnarray*}
\item[(c).] By the  Kuratowski-Ryll-Nardzewski selection theorem \cite[Theorem 18.13]{aliprantis06}, there exist measurable selectors for $x\mapsto \AA_1(x)$ and $x\mapsto\AA_2(x)$. In particular, $\mathbf{M}_1$ and $\mathbf{M}_2$ are nonempty.
\end{enumerate}
\end{remark}

\section{Occupation measures}\label{section3}
\setcounter{equation}{0}

\subsection{Occupation measures of the policies}

\paragraph{Definition of the occupation measures.}
Given an initial distribution $\eta\in\bcal{P}(\XX)$ and policies $(\pi_1,\pi_2)\in\mathbf{\Pi}_1\times\mathbf{\Pi}_2$ of the players, the occupation measure gives the total expected discounted time spent by the state-action process in a given subset of $\XX\times\AA$. 

\begin{definition}\label{def-occupation-measure}
Given an initial distribution $\eta\in\bcal{P}(\XX)$ and a pair of policies $(\pi_1,\pi_2)\in\mathbf{\Pi}_1\times\mathbf{\Pi}_2$, 
the associated \textit{occupation measure}   $\mu_{\eta,\pi_1,\pi_2}\in\bcal{P}(\XX\times\AA)$ is defined, for  $D\in\bfrak{X}$ and $(D_1,D_2)\in\mathfrak{B}(\AA_1)\times\mathfrak{B}(\AA_2)$ as
$$\mu_{\eta,\pi_1,\pi_2}(D\times D_1\times D_2)= (1-\beta)\sum_{t=0}^\infty \beta^t
\mathbb{P}_{\eta,\pi_1,\pi_2}\{X_t\in D, A_t\in D_1,B_t\in D_2\}.$$
The set of all occupation measures is denoted by $\bcal{O}_\eta=\{\mu_{\eta,\pi_1,\pi_2}:(\pi_1,\pi_2)\in\mathbf{\Pi}_1\times\mathbf{\Pi}_2\}\subseteq\bcal{P}(\XX\times\AA)$.
\end{definition}

As a direct consequence of Definition \ref{def-occupation-measure} and Assumption  \ref{Assumption-reward} we have that the total expected payoffs of the pair of policies $(\pi_1,\pi_2)\in\mathbf{\Pi}_1\times\mathbf{\Pi}_2$  for the initial distribution $\eta\in\bcal{P}(\XX)$ equal

\begin{equation}\label{eq-R-as-a-sum}
R_i(\eta,\pi_1,\pi_2)=\int_{\XX\times\AA} r_id\mu_{\eta,\pi_1,\pi_2}=
\int_{\XX\times\AA_1} r^1_id\mu^{\XX\times\AA_1}_{\eta,\pi_1,\pi_2}+\int_{\XX\times\AA_2} r^2_id\mu^{\XX\times\AA_2}_{\eta,\pi_1,\pi_2}
\end{equation}
for any $i=1,2$. Similar equalities hold for the constraints $C_i(\eta,\pi_1,\pi_2)$ of the players.

\paragraph{Occupation measures of stationary Markov policies.}
Given a pair of  stationary Markov policies $(\pi_1,\pi_2)\in\mathbf{M}_1\times\mathbf{M}_2$, 
we define the stochastic kernel $Q_{\pi_1,\pi_2}$ on $\XX$ given $\XX$ as
\begin{eqnarray}
Q_{\pi_1,\pi_2}(D|x)&=& \int_{\AA} Q(D|x,a_1,a_2)\pi_1(da_1|x)\pi_2(da_2|x) \label{eq-def-Qpi1pi2}
\end{eqnarray}
for $x\in\XX$ and $D\in\bfrak{X}$.
We denote by  $Q^t_{\pi_1,\pi_2}$ the $t$-th composition of $Q_{\pi_1,\pi_2}$ with itself
and we make the convention that $Q^0_{\pi_1,\pi_2}(\cdot|x)=\delta_x(\cdot)$.
It is then easily shown that the 
$\XX$-marginal of the occupation measure of $(\pi_1,\pi_2)\in\mathbf{M}_1\times\mathbf{M}_2$ for the initial distribution $\eta\in\bcal{P}(\XX)$ is\begin{equation}\label{Def-marginal-occup-measure}
\mu_{\eta,\pi_1,\pi_2}^\XX =(1-\beta)\sum_{t=0}^\infty \beta^t \eta Q^t_{\pi_1,\pi_2}
=\eta\Big[(1-\beta)\sum_{t=0}^\infty \beta^t  Q^t_{\pi_1,\pi_2}\Big].
\end{equation}
Therefore, we have
\begin{equation}
\mu_{\eta,\pi_{1},\pi_{2}}(dy,da_{1},da_{2})= \pi_{1}(da_{1}|y) \pi_{2}(da_{2}|y) \mu_{\eta,\pi_1,\pi_2}^\XX(dy),\label{eq-occupation-marginal-conditional}
\end{equation}
and so 
$\mu^{\XX\times\AA_1}_{\eta,\pi_{1},\pi_{2}}=\mu_{\eta,\pi_1,\pi_2}^\XX\otimes \pi_{1}$ and $
\mu^{\XX\times\AA_2}_{\eta,\pi_{1},\pi_{2}}=\mu_{\eta,\pi_1,\pi_2}^\XX\otimes \pi_{2}
$.

\paragraph{Basic results on occupation measures.}
Our next result summarizes the main results on the occupation measures of history-dependent policies. We note that statements (i) and (iii)--(iv) below are quite standard, while (ii) makes use of the ARAT nature of the game.

\begin{proposition}
\label{occup-meas-X-eqiv-lambda}
Let  $\eta\in\bcal{P}(\XX)$ be an initial distribution and let $(\pi_1,\pi_2)\in\mathbf{\Pi}_1\times\mathbf{\Pi}_2$ be an arbitrary pair of history-dependent policies.
\begin{enumerate}
\item[(i).] The occupation measure $\mu_{\eta,\pi_1,\pi_2}$ satisfies the linear equations (written in $\mu$) $$\mu(\KK)=1\quad\hbox{and}\quad
\mu^\XX=(1-\beta)\eta+\beta  \mu Q\qquad\hbox{for $\mu\in\bcal{P}(\XX\times\AA)$.}$$
\item[(ii).]  There exists a pair $(\pi'_1,\pi'_2)\in\mathbf{M}_1\times\mathbf{M}_2$ of stationary Markov policies  such that 
$$\mu^{\XX\times\AA_1}_{\eta,\pi_1,\pi_2}=\mu^{\XX\times\AA_1}_{\eta,\pi'_1,\pi'_2}\quad\hbox{and}\quad 
\mu^{\XX\times\AA_2}_{\eta,\pi_1,\pi_2}=\mu^{\XX\times\AA_2}_{\eta,\pi'_1,\pi'_2}.$$
If  $\pi_i\in\mathbf{M}_i$  then we can choose $\pi'_i=\pi_i$.
\item[(iii).]  If $\eta\in\bcal{P}_\lambda(\XX)$ (resp., $\eta\in\bcal{P}^e_\lambda(\XX)$) then $\mu^\XX_{\eta,\pi_1,\pi_2}\in\bcal{P}_\lambda(\XX)$ (resp., $\mu^\XX_{\eta,\pi_1,\pi_2}\in\bcal{P}^e_\lambda(\XX)$).
\end{enumerate}
\end{proposition}
\textbf{Proof.} (i). 
It is  well known in the theory of MDPs that the occupation measure satisfies the linear constraint in the statement of the proposition; see, for instance, \cite[Theorem 6.3.7]{hernandez96}.
\\[5pt]\noindent
(ii).
Given $i=1,2$, 
since the probability measure $\mu_{\eta,\pi_1,\pi_2}^{\XX\times\AA_i}$ is supported on the set $\KK_i$, by Lemma~\ref{lemma-disintegration} it follows that there exists $\pi'_i\in\mathbf{M}_i$ such that 
\begin{eqnarray}
\mu_{\eta,\pi_1,\pi_2}^{\XX\times\AA_i}(dx,da_i)=\mu_{\eta,\pi_1,\pi_2}^\XX (dx)\pi'_i(da_i|x)
\label{eq-tool-24bis}
\end{eqnarray}
(clearly, if  $\pi_i\in\mathbf{M}_i$ then we can let $\pi'_i=\pi_i$).
This implies that 
\begin{equation*}
\mu_{\eta,\pi_1,\pi_2}^\XX=(1-\beta)\eta +\beta \mu_{\eta,\pi_1,\pi_2} Q
=(1-\beta)\eta +\mu_{\eta,\pi_1,\pi_2}^{\XX}Q_{\pi'_1,\pi'_2}
\end{equation*}
recalling Assumption \ref{Assumption-transition-Q}.
However, the solution of the equation $\gamma=(1-\beta)\eta+\beta \gamma Q_{\pi'_1,\pi'_2}$ for $\gamma\in\bcal{P}(\XX)$ is unique
and   (recall \eqref{Def-marginal-occup-measure}) it equals $\mu_{\eta,\pi'_1,\pi'_2}^\XX.$
We have thus shown that $\mu_{\eta,\pi_1,\pi_2}^\XX=\mu^\XX_{\eta,\pi'_1,\pi'_2}$ and by  
\eqref{eq-occupation-marginal-conditional}--\eqref{eq-tool-24bis}
it follows that
$\mu_{\eta,\pi_1,\pi_2}^{\XX\times\AA_1}= \mu_{\eta,\pi'_1,\pi_2'}^{\XX\times\AA_1}$ and $
\mu_{\eta,\pi_1,\pi_2}^{\XX\times\AA_2}=\mu_{\eta,\pi_1',\pi_2'}^{\XX\times\AA_2}$, 
as we wanted to prove. 
\\[5pt]\noindent
(iii). The first statement easily follows because $\eta\ll\lambda$ and the transitions are absolutely continuous with respect to $\lambda$; recall Assumption \ref{Assumption-transition-Q}.
Regarding the second statement, observe from \eqref{Def-marginal-occup-measure} that
$(1-\beta) \eta\leq \mu^\XX_{\eta,\pi_1,\pi_2}$ and since, by hypothesis, we have $\eta\sim\lambda$,  we readily get the result.
\hfill$\Box$\\[10pt]
\indent
As a direct consequence of \eqref{eq-R-as-a-sum} and Proposition \ref{occup-meas-X-eqiv-lambda}(ii) we have the following fact.

\begin{corollary}\label{cor-Markov-sufficient}
For every initial distribution $\eta\in\bcal{P}(\XX)$ and any pair $(\pi_1,\pi_2)\in\mathbf{\Pi}_1\times\mathbf{\Pi}_2$, there exist stationary Markov policies 
$(\pi'_1,\pi'_2)\in\mathbf{M}_1\times\mathbf{M}_2$ such that for any $i=1,2$ we have
$$R_i(\eta,\pi_1,\pi_2)=R_i(\eta,\pi'_1,\pi'_2)\quad\hbox{and}\quad C_i(\eta,\pi_1,\pi_2)=C_i(\eta,\pi'_1,\pi'_2).$$
If $\pi_1\in\mathbf{M}_1$ (respectively, $\pi_2\in\mathbf{M}_2$) then the result holds for $\pi'_1=\pi_1$ (respectively, $\pi'_2=\pi_2$).
\end{corollary}
\begin{remark}\label{rem-Markov-sufficient}
By Corollary \ref{cor-Markov-sufficient},  the Slater condition  is equivalent to the  apparently weaker condition that for any $\pi_1\in\mathbf{M}_1$ there exists $\pi_2\in\mathbf{\Pi}_2$ with $C_2(\eta,\pi_1,\pi_2)>\rho_2$, and symmetrically for player 1.
\end{remark}

\paragraph{Stationary Markov policies and Young measures.}

In  Proposition \ref{occup-meas-X-eqiv-lambda}(iii) we have shown that, given any initial distribution $\eta\in\bcal{P}_\lambda(\XX)$, the state process
$\{X_t\}_{t\ge0}$ visits any $\lambda$-null set  with probability zero. Therefore, we can give an alternative definition of the set $\mathbf{M}_i$ of stationary Markov policies for player~$i$ 
by letting  $\mathbf{M}_i$ to be  the family of 
stochastic kernels $\pi_i$ on $\AA_i$ given~$\XX$ satisfying
\begin{equation}\label{eq-piAx}
\pi_i(\AA_{i}(x)|x)=1\quad\hbox{for $\lambda$-almost every $x\in\XX$}.
\end{equation} 
This definition does not alter any of the properties of Markov policies seen so far.
 
Let $L^{1}(\XX,\bfrak{X},\lambda)$ be the family of real-valued measurable functions which are $\lambda$-integrable (where, as usual, we identify functions which are equal $\lambda$-a.s.).  When endowed with the $\|\cdot\|_1$-norm, and as a consequence of Assumption \ref{Assumption-state}, we have that $L^{1}(\XX,\bfrak{X},\lambda)$ becomes a separable Banach space.

On the set  $\mathbf{M}_i$ we define the following equivalence relation. Given $\pi_i,\pi_i'\in\mathbf{M}_i$ we say that
$\pi_i\approx\pi_i'$ when $ \pi_i(\cdot|x)=\pi_i'(\cdot|x)$ for $ \lambda$-almost every $x\in\XX$.
Let $\bcal{Y}_i$ be the family of equivalence classes of this relation.
Each element of $\bcal{Y}_i$ is referred  to as a Young measure.
We equip the family of Young measures $\bcal{Y}_i$ with the narrow (stable) topology, which is the coarsest topology on $\bcal{Y}_i$ which makes the following  mappings continuous: $$\pi_i\mapsto \int_{\XX}\int_{\AA_i}f(x,a_i) \pi_i (da_i|x)\lambda(dx),$$
for any $f\in\car(\XX\times\AA_i,\RR)$ such that for some $\Phi$ in $L^1(\XX,\bfrak{X},\lambda)$
we have $|f(x,a_i)| \leq \Phi(x)$ for every $(x,a_{i})\in\XX\times\AA_i$; see  \cite[Theorem 2.2]{balder88}.
By \cite[Lemma 1]{balder91}, the set $\bcal{Y}_i$ endowed with the narrow topology  becomes a compact metric space. 

\begin{proposition}\label{proposition-uniqueness-Y}
Given any two pairs of stationary Markov policies $(\pi_1,\pi_2),(\pi'_1,\pi'_2)\in\mathbf{M}_1\times\mathbf{M}_2$ and an initial distribution $\eta\in\bcal{P}(\XX)$, we have the following results.
\begin{itemize}
\item[(i).] 
If  $\pi_i\approx\pi'_i$ (for  $i=1,2$) and $\eta\in\bcal{P}_\lambda(\XX)$, then $\mu_{\eta,\pi_1,\pi_2}=\mu_{\eta,\pi'_1,\pi'_2}$.
\item[(ii).] 
If  $\mu_{\eta,\pi_1,\pi_2}=\mu_{\eta,\pi'_1,\pi'_2}$ and $\eta\in\bcal{P}^e_\lambda(\XX)$, then $\pi_i\approx\pi'_i$ for each $i=1,2$.
\end{itemize}
\end{proposition}
\textbf{Proof.}
(i). Assume  that $\eta\ll\lambda$ and  $\pi_i\approx\pi'_i$ for $i=1,2$.
Observe that
$Q_{\pi_1,\pi_2}(\cdot | x)= Q_{\pi'_1,\pi'_2}(\cdot | x)$ for $ \lambda$-almost every $x\in\XX$, and that
$\eta Q^t_{\pi_1,\pi_2} \ll \lambda $ for every $t\ge0$. 
Then, it can be easily shown by induction  that $\eta Q^t_{\pi_1,\pi_2}=\eta Q^t_{\pi'_1,\pi'_2}$ for every $t\ge0$ implying ---by \eqref{Def-marginal-occup-measure}--- that
$\mu_{\eta,\pi_1,\pi_2}^\XX=\mu_{\eta,\pi'_1,\pi'_2}^\XX$. 
Once we know that the $\XX$-marginals of $\mu_{\eta,\pi_1,\pi_2}$ and $\mu_{\eta,\pi'_1,\pi'_2}$ coincide,
equality of the occupation measures follows from  \eqref{eq-occupation-marginal-conditional} because $\pi_i(\cdot|y)$ and $\pi_i'(\cdot|y)$ coincide for every $y\in\XX$ on a set of $\lambda$-probability one, which is also a set of $\mu^\XX_{\eta,\pi_1,\pi_2}=\mu^\XX_{\eta,\pi'_1,\pi_2'}$-probability one (recall Proposition \ref{occup-meas-X-eqiv-lambda}(iii)).
\\[5pt]\noindent
(ii).
Suppose now that $\eta\sim\lambda$ and $\mu_{\eta,\pi_1,\pi_2}=\mu_{\eta,\pi'_1,\pi'_2}$. We have
$$\mu^{\XX}_{\eta,\pi_1,\pi_2}(dy) \pi_i(da_i|y)= \mu^{\XX\times\AA_i}_{\eta,\pi'_1,\pi'_2}(dy,da_{i})=\mu^{\XX}_{\eta,\pi'_1,\pi'_2}(dy) \pi'_i(da_i|y)=\mu^{\XX}_{\eta,\pi_1,\pi_2}(dy) \pi'_i(da_i|y),$$
for $i=1,2$. Since disintegration of $\mu^{\XX\times\AA_i}_{\eta,\pi_1,\pi_2}$ is unique up to sets of $\mu^\XX_{\eta,\pi_1,\pi_2}$-probability zero and, hence, sets of $\lambda$-probability zero (recall Proposition \ref{occup-meas-X-eqiv-lambda}(iii)), it follows
that $\pi_i\approx\pi'_i$.
\hfill$\Box$\\[10pt]\indent
Proposition \ref{proposition-uniqueness-Y}(i) above shows that, whenever $\eta\in\bcal{P}_\lambda(\XX)$, stationary Markov policies in the same class of equivalence yield the same occupation measures. 
Therefore, in case that $\eta\in\bcal{P}_\lambda(\XX)$, we will henceforth refer to the sets of Young measures $\bcal{Y}_1$ and $\bcal{Y}_2$ as to the stationary Markov policies of the players, with a slight abuse of terminology.

\subsection{Some continuity facts}

Let  $L^{\infty}(\XX,\bfrak{B}(\XX),\lambda)$ be the set of $\lambda$-essentially bounded measurable real-valued functions on $\XX$ endowed with the weak$^{*}$ topology (we identify  functions which are equal $\lambda$-a.s.). Let $\|v\|$ be  the essential supremum of $v\in L^{\infty}(\XX,\bfrak{B}(\XX),\lambda)$. The next  result follows from Assumption \ref{Assumption-transition-Q}.

\begin{lemma}\label{lem-v-Qv}
If $v\in L^{\infty}(\XX,\bfrak{B}(\XX),\lambda)$ then $Qv\in\car_b(\XX\times\AA,\RR)$.
\end{lemma}

\begin{lemma}
\label{Nowak2019}
Suppose that 
 $(\pi_{1,n},\pi_{2,n})\rightarrow (\pi_1,\pi_2)$ in $\bcal{Y}_1\times\bcal{Y}_2$ and $v_n\wstar v$  in $L^{\infty}(\XX,\bfrak{X},\lambda)$ as $n\rightarrow\infty$. Under these conditions, for  any $t\ge0$ we have 
\begin{align}
Q^t_{\pi_{1,n},\pi_{2,n}} v_{n} \wstar Q^t_{\pi_{1},\pi_{2}}v\quad\hbox{as $n\rightarrow\infty$ in $L^{\infty}(\XX,\bfrak{X},\lambda)$}.
\end{align}
\end{lemma}
\textbf{Proof.} 
The result is trivial for $t=0$. The case $t=1$ can be shown by using similar arguments as in the proof of Lemma 4.1 in \cite{nowak19}.
Once we know that $Q_{\pi_{1,n},\pi_{2,n}} v_{n} \wstar Q_{\pi_{1},\pi_{2}}v$,  the stated result for any integer $t\ge2$ follows easily.
\hfill$\Box$
\begin{proposition}
\label{Convergence-mu}
The mappings from $\bcal{P}_{\lambda}(\XX)\times\bcal{Y}_{1}\times\bcal{Y}_{2}$ to $\bcal{P}(\XX\times\AA_i)$defined by $(\eta,\pi_1,\pi_2)\mapsto \mu_{\eta,\pi_1,\pi_2}^{\XX\times\AA_i}$ are continuous for  $i=1,2$
where $\bcal{P}_\lambda(\XX)$ is endowed with the metric of total variation and $\bcal{P}(\XX\times\AA_i)$ with the $ws$-topology.
As a consequence, the $\RR$- and $\RR^p$-valued mappings 
$$(\eta,\pi_1,\pi_2)\mapsto R_i(\eta,\pi_1,\pi_2)\quad\hbox{and}\quad (\eta,\pi_1,\pi_2)\mapsto C_i(\eta,\pi_1,\pi_2)$$ are also continuous on
$\bcal{P}_{\lambda}(\XX)\times\bcal{Y}_{1}\times\bcal{Y}_{2}$ for $i=1,2$.
\end{proposition}
\textbf{Proof:}
Fix $i=1,2$.  Since $\bcal{P}_{\lambda}(\XX)\times\bcal{Y}_1\times\bcal{Y}_2$ is a metric space, we will check continuity  by proving sequential continuity.
Consider a convergent sequence  $(\eta_n,\pi_{1,n},\pi_{2,n})\rightarrow (\eta,\pi_1,\pi_2)$ in  the product space $\bcal{P}_\lambda(\XX)\times\bcal{Y}_1\times\bcal{Y}_2$
and  $f\in\car_b(\XX\times\AA_{i},\RR)$.
Let us denote by $v_{n}$  and $v$ the bounded measurable functions on $\XX$ defined by 
$v_{n}(x)=\int_{\AA_{i}} f(x,a_i) \pi_{i,n}(da_i|x)$ and $v(x)=\int_{\AA_{i}} f(x,a_i) \pi_{i}(da_i|x)$. Since
$\pi_{i,n}\rightarrow \pi_{i}$ in $\bcal{Y}_{i}$ then
$v_n \wstar v$ in $L^{\infty}(\XX,\bfrak{X},\lambda)$.
We can now apply Lemma \ref{Nowak2019}  to conclude that for any $t\ge0$ we have
$Q^t_{\pi_{1,n},\pi_{2,n}} v_{n} \wstar Q^t_{\pi_{1},\pi_{2}} v$ in $L^{\infty}(\XX,\bfrak{B}(\XX),\lambda)$.
Moreover, by hypothesis $\eta_{n}\rightarrow\eta$ in total variation, and so it is easily seen that  ${d\eta_{n}}/{d\lambda} \rightarrow {d\eta}/{d\lambda}$ in
$L^{1}(\XX,\bfrak{X},\lambda)$.
Therefore, 
\begin{eqnarray*}
\lim_{n\rightarrow\infty} \int_\XX Q_{\pi_{1,n},\pi_{2,n}}^{t}v_n(x)\eta_{n}(dx) & =&  \lim_{n\rightarrow\infty} \int_\XX Q_{\pi_{1,n},\pi_{2,n}}^{t}v_n(x)
\frac{d\eta_{n}}{d\lambda}(x) \lambda(dx) \\
& = &  \int_\XX Q_{\pi_{1},\pi_{2}}^{t}v(x) \frac{d\eta}{d\lambda}(x) \lambda(dx)= \int_{\XX} Q_{\pi_{1},\pi_{2}}^{t} v(x) \eta(dx),
\end{eqnarray*}
by using Proposition 3.13(iv) in \cite{brezis11}.
Therefore, since $f$ is bounded, by dominated convergence
\begin{align*}
\lim_{n\rightarrow\infty} \int_{\XX\times\AA_{i}} & f(x,a_i) \mu_{\eta_{n}\pi_{1,n},\pi_{2,n}}^{\XX\times\AA_{i}} (dx,da_i)
= (1-\beta)  \sum_{t=0}^{\infty} \beta^{t} \lim_{n\rightarrow\infty}  \int_\XX  Q_{\pi_{1,n},\pi_{2,n}}^{t}v_{n}(x)\eta_{n}(dx) \\
& = (1-\beta)   \sum_{t=0}^{\infty} \beta^{t}  \int_\XX Q_{\pi_{1},\pi_{2}}^{t}v(x)\eta(dx) =  \int_{\XX\times\AA_{i}} f(x,a_i) \mu_{\eta,\pi_{1},\pi_{2}}^{\XX\times\AA_{i}} (dx,da_i),
\end{align*}
which establishes 
the continuity of $(\eta,\pi_1,\pi_2)\mapsto \mu_{\eta,\pi_1,\pi_2}^{\XX\times\AA_i}$. Continuity of the reward and constraint functions  follows because they are defined by integration of bounded Carath\'eodory functions. 
\hfill $\Box$

\subsection{The sets $\mathcal{D}_{\eta}$ and $\mathcal{D}_{\eta,i}$}

Based on the result of Proposition \ref{occup-meas-X-eqiv-lambda}(i), given an initial distribution $\eta\in\bcal{P}(\XX)$ we define 
$$
\mathcal{D}_{\eta}=  \big\{\mu\in\bcal{P}(\XX\times\AA) : \mu(\KK)=1 \text{ and } \mu^{\XX}=(1-\beta)\eta+\beta \mu Q \big\}.
$$
The set of occupation measures satisfies $\bcal{O}_\eta\subseteq\mathcal{D}_\eta$.
Note that any $\mu\in\mathcal{D}_\eta$ can be disintegrated  as 
\begin{equation}\label{eq-disintegrate-linear}
\mu(dx,da_1,da_2)=\pi(da_1,da_2|x)\mu^\XX(dx)
\end{equation}
for some stochastic kernel $\pi$ on $\AA$ given $\XX$ which satisfies $\pi(\AA(x)|x)=1$ for all $x\in\XX$, although this $\pi$ might not correspond to a pair $(\pi_1,\pi_2)\in\mathbf{M}_1\times\mathbf{M}_2$ of stationary Markov policies of the players. The stochastic kernel  $Q_\pi$ on $\XX$ given $\XX$ is then defined as $Q_\pi(dy|x)=\int_{\AA} Q(dy|x,a)\pi(da|x)$ 
for $x\in\XX$,
and $Q^t_\pi$ for $t\ge0$ denotes the $t$-th composition of $Q_{\pi}$, with $Q^0_\pi(dy|x)=\delta_x(dy)$.  Clearly, 
\begin{equation}\label{eq-marginal-linear-equation}
\mu^\XX=(1-\beta)\sum_{t=0}^\infty \beta^t \eta Q^t_\pi.
\end{equation}
For $\eta\in\bcal{P}(\XX)$, we define $\mathcal{D}_{\eta,i}$ as the set of  $(\XX\times\AA_i)$-marginals of the measures in $\mathcal{D}_{\eta}$ for $i=1,2$: 
$$\mathcal{D}_{\eta,i}=\{\mu^{\XX\times\AA_i}:\mu\in\mathcal{D}_{\eta}\}\subseteq\bcal{P}(\XX\times\AA_i).$$

In our next result we establish compactness of the above defined sets.  Such result is known in the literature when
$\mathbf{X}$ is a Borel space  whereas, in this paper,  $\mathbf{X}$ is an abstract measurable space.
\begin{proposition}
\label{Compactness-set-D}
Given an initial distribution $\eta\in\bcal{P}(\XX)$, the sets $\mathcal{D}_{\eta}$ and $\mathcal{D}_{\eta,i}$  are convex compact metric spaces when endowed with their respective $ws$-topologies.
\end{proposition}
\textbf{Proof.} 
Let us  first show that $\mathcal{D}_{\eta}$  is relatively compact in $\bcal{P}(\XX\times\AA)$.
To this end, consider a decreasing sequence $\{h_{j}\}_{j\in\NN}$ of functions in $\car_b(\XX\times\AA,\RR)$ such that 
$h_{j}(x,a_1,a_2)\downarrow0$ for any $(x,a_1,a_2)\in\XX\times\AA$.
By Proposition \ref{Criterion-relatively-ws-compactness},  to prove relative compactness we need to show that
\begin{align}\label{eq-relative-compactness-toprove}
\lim_{ j \rightarrow\infty} \sup_{\mu\in\mathcal{D}_{\eta}} \int_{\XX\times\AA} h_{j}d\mu \ =0 .
\end{align}
According to \eqref{eq-disintegrate-linear} and \eqref{eq-marginal-linear-equation}, for any $\mu\in\mathcal{D}_{\eta}$ there exists a stochastic kernel $\pi_{\mu}$ on $\AA$ given $\XX$ satisfying $\pi_{\mu}(\AA(x)|x)=1$ for all $x\in\XX$ such that 
\begin{align}
 \int_{\KK} h_{j}(x,a) \mu(dx,da) \leq \int_{\XX} g_{j}^{\mu}(x) \eta(dx)\quad\hbox{where}\quad
 g_{j}^{\mu}(x)= (1-\beta)\sum_{t=0}^{\infty} \beta^{t} Q_{\pi_{\mu}}^{t} f_{j}(x)
 \label{Compactness-set-D-temp-1}
\end{align} 
with $f_{j}(x)=\sup_{a\in \AA} h_{j}(x,a)$.
Combining Lemma 10.1 and Theorem 12.1 in \cite{schal75b}, we obtain that 
$\{f_{j}\}_{j\in\NN}$ is a decreasing sequence  of bounded measurable functions defined on $\XX$ which satisfies 
$\ds \lim_{ j \rightarrow\infty} f_{j}(x)=\sup_{a\in\AA} \lim_{j\rightarrow\infty} h_{j}(x,a)=0 $ for any $x\in\XX$.
Define the functions $\bar{f}_{j,t}$ for $j,t\in\NN$ on $\XX$ by 
$$\bar{f}_{j,t}(x)=\sup_{a\in\AA} Q\bar{f}_{j,t-1}(x,a)\quad\hbox{and}\quad
\bar{f}_{j,0}(x)=f_{j}(x)$$
 for $x\in\XX$, $t\geq 1$,  and $j\in\NN$.
Observe that the functions $\{\bar{f}_{j,t}\}_{j,t\in\NN}$ do not depend on $\mu\in\mathcal{D}_\eta$.
Clearly, we have 
$(Q_{\pi_{\mu}})^{t}f_{j}(x)\leq \bar{f}_{j,t}(x) $  for any $j,t\in\NN$. From \eqref{Compactness-set-D-temp-1}, this implies for any $x\in\XX$ and $j\in\NN$ that
$g_{j}^{\mu}(x) \leq (1-\beta)\sum_{t=0}^{\infty} \beta^{t}  \bar{f}_{j,t}(x)$.
It follows that
\begin{align}
\sup_{\mu\in\mathcal{D}} \int_{\XX\times\AA} h_{j}(x,a) \mu(dx,da) \leq (1-\beta)\sum_{k=0}^{\infty} \beta^{k}  \int_{\XX} \bar{f}_{j,k}(x) \eta(dx).
\label{Compactness-set-D-temp-3}
\end{align}
From Lemma 10.1 and Theorem 12.1 in \cite{schal75b}, and recalling Assumption \ref{Assumption-transition-Q}, it can be shown  by induction that for any $k\in\NN$, $\{\bar{f}_{j,k}\}_{j\in\NN}$ is a decreasing sequence  of bounded measurable functions satisfying
$\lim_{ j \rightarrow\infty} \bar{f}_{j,k}(x)=0$ for any $x\in\XX$  and $k\in\NN$.
 Therefore, taking the limit in \eqref{Compactness-set-D-temp-3} and by using the monotone convergence theorem, we get \eqref{eq-relative-compactness-toprove}
 and thus $\mathcal{D}_{\eta}$ is indeed relatively compact. 

Once we know that $\mathcal{D}_{\eta}$ is relatively compact we use \cite[Theorem 5.2]{balder01} and the fact that $\AA$ is a compact metric space to conclude that $\mathcal{D}^\XX_{\eta}=\{\mu^{\XX} : \mu\in \mathcal{D}_{\eta}\}$ is relatively compact in the $s$-topology of $\bcal{P}(\XX)$. Since $\bfrak{X}$ is countably generated, we deduce from  \cite[Proposition 2.3]{balder01} that $\mathcal{D}_{\eta}$ is metrizable. 

To conclude the proof of compactness, it remains to show that $\mathcal{D}_{\eta}$ is closed.
To this end,  let $\{\mu_{n}\}_{n\in\NN}$ be a sequence in $\mathcal{D}_{\eta}$ that converges in the $ws$-topology to
 $\mu\in\bcal{P}(\XX\times\AA)$. For any bounded measurable function $f:\XX\rightarrow\RR$ we have
$$
\int_\XX fd\mu^\XX_n=\int_{\XX\times\AA} fd\mu_n=(1-\beta)\int_\XX fd\eta +\beta \int_{\XX\times\AA}Qf(x,a_1,a_2) \mu_n
(dx,da_1,da_2),$$
where we still use the notation $f$ for the function $(x,a_1,a_2)\mapsto f(x)$ which is in $\car_b(\XX\times\AA,\RR)$.
By Lemma \ref{lem-v-Qv}, we can take the limit as $n\rightarrow\infty$ in order to obtain
$$ \int_\XX fd\mu^\XX=(1-\beta)\int_\XX fd\eta +\beta \int_{\XX\times\AA} Qf(x,a_1,a_2) \mu (dx,da_1,da_2)$$
and, in particular, $\mu^{\XX}=(1-\beta)\eta+\beta \mu Q$. To show that $\mu\in\mathcal{D}_{\eta}$ it remains to prove that $\mu(\KK)=1$. 
The mapping $(x,a_1,a_2)\mapsto -\mathbf{I}_{\KK}(x,a_1,a_2)$ is measurable on $\XX\times\AA$ and it is such that
$(a_1,a_2)\mapsto -\mathbf{I}_{\KK}(x,a_1,a_2)=-\mathbf{I}_{\AA(x)}(a_1,a_2)$
is lower semicontinuous for any fixed $x\in\XX$ because $\AA(x)$ is compact. Thus,
$-\mathbf{I}_{\KK}$ is a normal integrand \cite[p. 502]{balder01} and so by \cite[Theorem 3.1.(c)]{balder01} we have  
$\limsup_{n} \mu_{n}(\KK)\leq \mu(\KK)$,  implying that $\mu(\KK)=1$. 
This concludes the proof that $\mu\in\mathcal{D}_{\eta}$.  Finally, observe that convexity of $\mathcal{D}_{\eta}$ is a straightforward consequence of its definition.

 It is easy to check that the mapping from $\bcal{P}(\XX\times\AA)$ to $\bcal{P}(\XX\times\AA_{i})$ that associates to each $\mu\in\bcal{P}(\XX\times\AA)$ its marginal probability measure $\mu^{\XX\times\AA_i}$ is continuous for the $ws$-topologies. Hence~$\mathcal{D}_{\eta,i}$ is compact in $\bcal{P}(\XX\times\AA_{i})$. Recalling that $\bfrak{X}$ is countably generated and noting that the set of $\XX$-marginal probability measures of $\mathcal{D}_{\eta,i}$ is precisely $\mathcal{D}^\XX_{\eta}$, which has already been shown to be relatively compact for the $s$-topology, we get from   \cite[Proposition 2.3]{balder01} that
$\mathcal{D}_{\eta,i}$ is metrizable.

Finally, regarding convexity, 
let $\hat{\gamma}$ and $\bar{\gamma}$ in $\mathcal{D}_{\eta,i}$ and fix $0\le \alpha\le1$. There exist $\hat{\mu}$ and $\bar{\mu}$ in $\mathcal{D}_{\eta}$ such that
$ \hat{\gamma}=\hat{\mu}^{\XX\times\AA_{i}}$ and $ \bar{\gamma}=\bar{\mu}^{\XX\times\AA_{i}}$. By convexity of $\mathcal{D}_{\eta}$
we have  $\mu=\alpha\hat{\mu}+(1-\alpha)\bar{\mu}\in\mathcal{D}_{\eta}$ and
so $\alpha\hat{\gamma}+(1-\alpha)\bar{\gamma}=\mu^{\XX\times\AA_{i}}$ is in $\mathcal{D}_{\eta,i}$, showing convexity.
\hfill$\Box$

\begin{proposition}
\label{set-D-policies}
Consider an initial distribution $\eta\in\bcal{P}^{e}_{\lambda}(\XX)$.
\begin{itemize}
\item[(i).]
For any $\mu\in \mathcal{D}_{\eta}$ then there exists a unique $(\pi_{1},\pi_{2})\in\bcal{Y}_{1}\times\bcal{Y}_{2}$ satisfying
\begin{align}\label{eq-toprove-d}
\mu^{\XX\times\AA_{1}}= \mu_{\eta,\pi_{1},\pi_{2}}^{\XX\times\AA_{1}}\quad\hbox{and}\quad   \mu^{\XX\times\AA_{2}}= \mu_{\eta,\pi_{1},\pi_{2}}^{\XX\times\AA_{2}}.
\end{align}
\item[(ii).] We have $\mathcal{D}_{\eta,i}=\{\mu^{\XX\times\AA_i}_{\eta,\pi_1,\pi_2}:(\pi_1,\pi_2)\in\bcal{Y}_1\times\bcal{Y}_2\}$ for each $i=1,2$.
\item[(iii).] If $(\pi_1,\pi_2)$ and $(\pi'_1,\pi'_2)$ in $\bcal{Y}_1\times\bcal{Y}_2$ are such that 
$\mu^{\XX\times\AA_1}_{\eta,\pi_1,\pi_2}=\mu^{\XX\times\AA_1}_{\eta,\pi'_1,\pi'_2}$
then $\pi_1=\pi'_1$.  A symmetric result holds for the  marginal probability measures on $\XX\times\AA_2$.
\end{itemize}
\end{proposition}
\textbf{Proof.} (i). 
Since $\mu(\KK)=1$ we  have  $\mu^{\XX\times\AA_1}(\KK_1)=1$ and $\mu^{\XX\times\AA_2}(\KK_2)=1$. By Lemma~\ref{lemma-disintegration} there exist 
 $(\pi_{1},\pi_{2})\in\mathbf{M}_{1}\times\mathbf{M}_{2}$ such that 
$\mu^{\XX\times\AA_1}=\mu^{\XX}\otimes \pi_1$ and
$\mu^{\XX\times\AA_2}=\mu^{\XX}\otimes \pi_2$.
Now,  $\mu$ satisfies 
$$\mu^\XX=(1-\beta)\eta+\beta \mu Q=(1-\beta)\eta+\beta \mu^\XX Q_{\pi_1,\pi_2}$$ and so
 $\mu^\XX=\mu^\XX_{\eta,\pi_1,\pi_2}$.
The above stochastic kernels $\pi_{i}$ on $\AA_{i}$ given $\XX$  are unique $\mu^\XX$-almost surely. Since $\mu^\XX=\mu^\XX_{\pi_1,\pi_2}\sim\lambda$ (Proposition~\ref{occup-meas-X-eqiv-lambda}(iii)), it follows that uniqueness is $\lambda$-almost surely. This shows that there is indeed a unique pair 
$(\pi_1,\pi_2)\in \bcal{Y}_1\times\bcal{Y}_2$ with the above mentioned properties.
\\[5pt]\noindent
(ii). This result directly follows from part (i) and the definition of $\mathcal{D}_{\eta,i}$. 
\\[5pt]\noindent
(iii). By using arguments similar to those in (i) above, the result readily follows.
\hfill $\Box$

\section{Main results}\label{section4}
\setcounter{equation}{0}

In this section we address our results on the existence of Nash equilibria. First of all, in Section \ref{section41}, we will prove the existence of such equilibria for a game model with initial distribution in $\bcal{P}^e_\lambda(\XX)$, while in  Section \ref{section42} we will treat the general case a  game model $\mathcal{G}(\nu,\theta)$ satisfying Assumption \ref{Assump}.

\subsection{The case of an initial distribution $\eta\in\bcal{P}^e_\lambda(\XX)$}\label{section41}

Suppose that the initial distribution of the system $\eta\in\bcal{P}^e_\lambda(\XX)$ and that  the constraint constants $\rho=(\rho_1,\rho_2)\in\RR^p\times\RR^p$ that satisfy the Slater condition. 
In order to prove the existence of a constrained Nash equilibrium for $\mathcal{G}(\eta,\rho)$ we will suitably define a correspondence 
$\mathcal{D}_{\eta,1}\times\mathcal{D}_{\eta,2}\tto \mathcal{D}_{\eta,1}\times\mathcal{D}_{\eta,2}$ which will be shown to have a fixed point, from which we  will derive equilibrium stationary Markov policies. In fact, we shall construct this correspondence in two steps: as the composition
$\mathcal{H}_{\eta,\rho}\circ\mathcal{J}_{\eta}$ of a function
$\mathcal{J}_{\eta}:\mathcal{D}_{\eta,1}\times\mathcal{D}_{\eta,2}\rightarrow \bcal{Y}_1\times\bcal{Y}_2$ and a correspondence
$\mathcal{H}_{\eta,\rho}: \bcal{Y}_1\times\bcal{Y}_2\tto\mathcal{D}_{\eta,1}\times\mathcal{D}_{\eta,2}$.

\paragraph{The function $\mathcal{J}_{\eta}$.}
Given a player $i=1,2$  we define the function $\mathcal{J}_{\eta,i}:\mathcal{D}_{\eta,i}\to\bcal{Y}_{i}$ as follows. For any $\gamma\in\mathcal{D}_{\eta,i}$, let $\pi_i=\mathcal{J}_{\eta,i}(\gamma)$ be the unique 
$\pi_i\in\bcal{Y}_i$ satisfying $\gamma=\gamma^\XX\otimes \pi_i$. The existence and uniqueness of such decomposition is guaranteed by (ii) and (iii) in Proposition \ref{set-D-policies}. 
\begin{proposition}
\label{Correspondence-J}
Given $\eta\in\bcal{P}^{e}_{\lambda}(\XX)$, the functions $\mathcal{J}_{\eta,i}:\mathcal{D}_{\eta,i}\to\bcal{Y}_{i}$ are continuous for $i=1,2$.
\end{proposition}
\textbf{Proof:} 
For simplicity in the notation, we prove the case $i=1$. 
Consider a sequence $\{\gamma_{n}\}_{n\in\NN}$ in~$\mathcal{D}_{\eta,1}$ such that $\gamma_{n}\rightarrow\gamma$ for some $\gamma\in\mathcal{D}_{\eta,1}$.
Using Proposition \ref{set-D-policies}(ii) it follows that 
for each $n\in\NN$  there exist $(\pi_{1,n},\pi_{2,n})\in \bcal{Y}_{1}\times \bcal{Y}_{2}$ satisfying
$\gamma_{n} = \mu^{\XX\times\AA_{1}}_{\eta,\pi_{1,n},\pi_{2,n}}$ and there exist $(\pi_{1},\pi_{2})\in \bcal{Y}_{1}\times\bcal{Y}_{2}$ such that
$\gamma= \mu^{\XX\times\AA_{1}}_{\eta,\pi_{1},\pi_{2}}$.
We have, by definition, $\pi_{1,n}=\mathcal{J}_{\eta,1}(\gamma_n)$ and $\pi_1=\mathcal{J}_{\eta,1}(\gamma)$.
Our goal is to show that $\pi_{1,n}\rightarrow\pi_1$.
Since~$\bcal{Y}_1$ is a compact metric space, to prove the result it suffices to consider an arbitrary convergent subsequence $\{\pi_{1,n'}\}$ of $\{\pi_{1,n}\}$ with $\pi_{1,n'}\rightarrow\pi_1'\in\bcal{Y}_1$ and to show that, necessarily,  $\pi'_1=\pi_1$.  But~$\bcal{Y}_2$ being also a compact metric space, there exists a further subsequence ---which, without loss of generality, we will denote also by $n'$--- such that $\pi_{2,n'}\rightarrow\pi'_2$ for some $\pi'_2\in\bcal{Y}_2$.
By Proposition \ref{Convergence-mu} we have
$$\gamma_{n'}=\mu^{\XX\times\AA_{1}}_{\eta,\pi_{1,n'},\pi_{2,n'}}\rightarrow \mu^{\XX\times\AA_{1}}_{\eta,\pi'_1,\pi'_{2}}=\gamma.$$
Recalling Proposition \ref{set-D-policies}(iii),
this implies that $\pi_1=\pi'_1=\mathcal{J}_{\eta,1}(\gamma)$, and the proof is complete.
\hfill $\Box$
\\[10pt]
\indent
Based on this result, we can now write $\mathcal{J}_{\eta}=(\mathcal{J}_{\eta,1},\mathcal{J}_{\eta,2})$ which is therefore a continuous function from
$\mathcal{D}_{\eta,1}\times\mathcal{D}_{\eta,2}$ to $\bcal{Y}_1\times\bcal{Y}_2$ for any $\eta\in\bcal{P}^{e}_{\lambda}(\XX)$.

\paragraph{The correspondence $\mathcal{H}_{\eta,\rho}$.}

Recall that we are considering $\eta\in\bcal{P}^{e}_{\lambda}(\XX)$. Given any $\pi_2\in\bcal{Y}_2$ let 
$$
\mathcal{L}_{\eta,1}(\pi_{2}) =  \big\{
 \mu_{\eta,\pi_{1},\pi_{2}}^{\XX\times\AA_{1}} : \pi_{1}\in\bcal{Y}_{1} \big\}
\subseteq\mathcal{D}_{\eta,1}\subseteq\bcal{P}(\XX\times\AA_1),
$$
which is the set of $(\XX\times\AA_1)$-marginals of the occupation measures for the initial distribution $\eta$ and pairs of stationary Markov policies $(\pi_1,\pi_2)\in\bcal{Y}_1\times\bcal{Y}_2$ when  the policy $\pi_2\in\bcal{Y}_2$ of player 2 remains fixed and the policy $\pi_1\in\bcal{Y}_1$ of player 1 varies. Similarly, for any 
 $\pi\in\bcal{Y}_1$ we define
$$\mathcal{L}_{\eta,2}(\pi_{1}) =  \big\{ \mu_{\eta,\pi_{1},\pi_{2}}^{\XX\times\AA_{2}} : \pi_{2}\in\bcal{Y}_2\}\subseteq \mathcal{D}_{\eta,2}\subseteq\bcal{P}(\XX\times\AA_2).
$$
\begin{proposition}\label{prop-L1-convex-compact}
Given an initial distribution $\eta\in\bcal{P}^{e}_{\lambda}(\XX)$ and any $(\pi_1,\pi_2)\in\bcal{Y}_1\times\bcal{Y}_2$, the sets $\mathcal{L}_{\eta,1}(\pi_{2})$ and $\mathcal{L}_{\eta,2}(\pi_{1})$ are convex and compact when endowed with their respective $ws$-topologies.
\end{proposition}
\textbf{Proof.} We only make the proof for $\mathcal{L}_{\eta,1}(\pi_{2})$.
To prove convexity, we fix
two measures $\tilde{\gamma},\bar{\gamma}$ in $\mathcal{L}_{\eta,1}(\pi_2)$ and some $0\le\alpha\le1$. We want to prove that $\alpha\tilde\gamma+(1-\alpha)\bar\gamma\in\mathcal{L}_{\eta,1}(\pi_2)$. 
By definition of $\mathcal{L}_{\eta,1}(\pi_2)$, there exist $\tilde{\pi}_{1},\bar{\pi}_{1}\in\bcal{Y}_{1}$
such that 
$\tilde{\gamma}= \mu_{\eta,\tilde{\pi}_{1},\pi_{2}}^{\XX\times\AA_{1}}$ and $\bar{\gamma}= \mu_{\eta,\bar{\pi}_{1},\pi_{2}}^{\XX\times\AA_{1}}$.
By convexity of $\mathcal{D}_\eta$  (recall Proposition~\ref{Compactness-set-D}), we have that  $\mu$ defined by
$\mu=\alpha \mu_{\eta,\tilde\pi_1,\pi_2}+(1-\alpha)\mu_{\eta,\bar\pi_1,\pi_2}$ lies in
$\mathcal{D}_\eta$.
Therefore, by Proposition \ref{set-D-policies}(i), there exists a unique $(\pi^*_1,\pi^*_2)\in\bcal{Y}_1\times\bcal{Y}_2$ such that 
$\mu^{\XX\times\AA_1}=\mu^{\XX\times\AA_1}_{\eta,\pi^*_1,\pi^*_2}$ and $\mu^{\XX\times\AA_2}=\mu^{\XX\times\AA_2}_{\eta,\pi^*_1,\pi^*_2}$.
Clearly, we have
\begin{eqnarray}
\label{eq-tool-14}
\alpha\tilde\gamma+(1-\alpha)\bar\gamma= \mu^{\XX\times\AA_1}_{\eta,\pi^*_1,\pi^*_2}.
\end{eqnarray}
To get convexity, we need to show that $\pi^*_2=\pi_2$. Observe that 
\begin{eqnarray}
\alpha\tilde\gamma+(1-\alpha)\bar\gamma &=&  \big[\alpha \mu_{\eta,\tilde\pi_1,\pi_2}^{\XX}+(1-\alpha)\mu_{\eta,\bar\pi_1,\pi_2}^{\XX}\big] \otimes \pi_{2}
=\mu^{\XX}\otimes \pi_{2}=\mu^{\XX}_{\eta,\pi^*_1,\pi^*_2}\otimes \pi_{2} \label{eq-tool-14bis}\\
\mu^{\XX\times\AA_1}_{\eta,\pi^*_1,\pi^*_2}&=&\mu^{\XX}_{\eta,\pi^*_1,\pi^*_2}\otimes \pi^{*}_{2}.\label{eq-tool-15}
\end{eqnarray}
However, $\mu^{\XX}_{\eta,\pi^*_1,\pi^*_2}\sim\lambda$ from Proposition \ref{occup-meas-X-eqiv-lambda}(iii) and so, by uniqueness of the disintegration, we deduce from \eqref{eq-tool-14}--\eqref{eq-tool-15} 
that $\pi_2=\pi^*_2$. Hence, $\mathcal{L}_{\eta,1}(\pi_2)$ is convex.

Since $\mathcal{L}_{1,\eta}(\pi_2)\subseteq\mathcal{D}_{1,\eta}$, which is a compact metric space, in order to prove compactness of $\mathcal{L}_{\eta,1}(\pi_2)$ is suffices to show that it is closed. To this end, consider a sequence $\{\gamma_n\}$ in $\mathcal{L}_{\eta,1}(\pi_2)$ converging to some $\gamma\in\mathcal{D}_{\eta,1}$. We have 
$\gamma_n=\mu^{\XX\times\AA_1}_{\eta,\pi_{1,n},\pi_2}$
for some $\pi_{1,n}$. For some subsequence $n'$ and some $\pi_1\in\bcal{Y}_1$ we have $\pi_{1,n'}\rightarrow \pi_1$ (recall that $\bcal{Y}_1$ is compact) and, by Proposition \ref{Convergence-mu}, this implies 
$$\gamma_{n'}=\mu^{\XX\times\AA_1}_{\eta,\pi_{1,n'},\pi_2}\rightarrow \mu^{\XX\times\AA_1}_{\eta,\pi_{1},\pi_2}=\gamma.$$
This shows that $\gamma\in\mathcal{L}_{\eta,1}(\pi_2)$ and the proof that $\mathcal{L}_{\eta,1}(\pi_2)$ is compact is complete.\hfill$\Box$\\[10pt]\indent
%
We explain how the payoff of the players for the policies $(\pi_1,\pi_2)\in\bcal{Y}_1\times\bcal{Y}_2$ relates to  $\mathcal{L}_1(\pi_2)$.

\begin{remark}\label{remark-cost}
Recalling \eqref{eq-R-as-a-sum} observe that for an initial distribution $\eta\in\bcal{P}^{e}_{\lambda}(\XX)$,
when player 2 fixes a policy $\pi_2\in\bcal{Y}_2$,
the family of total expected discounted payoffs $R_1(\eta,\cdot,\pi_2)$ when player 1 uses stationary Markov policies in $\bcal{Y}_1$ is given by 
\begin{equation}\label{eq-payoff-gamma}
\int_{\XX\times\AA_1} r^1_1d\gamma +\int_{\XX\times\AA_2} r_1^2d(\gamma^\XX\otimes\pi_2)
\end{equation}
when $\gamma$ varies in $\mathcal{L}_{\eta,1}(\pi_2)$. The important fact is that~\eqref{eq-payoff-gamma} is linear in $\gamma$. The same result holds for the
constraint $C_1(\eta,\cdot,\pi_2)$ and, symmetrically, for the payoffs $R_2(\eta,\pi_1,\cdot)$ and $C_2(\eta,\pi_1,\cdot)$ of player~2 when player 1 fixes his policy $\pi_1\in\bcal{Y}_1$.
\end{remark}

Based on this remark, we define the following sets. Given the initial distribution  $\eta\in\bcal{P}^{e}_{\lambda}(\XX)$, the constraint constant $\rho_{1}\in\RR^{p}$, and a stationary Markov policy  $\pi_2\in\bcal{Y}_2$,
let $\mathcal{A}_{\eta,\rho_1,1}(\pi_2)$ be the set of $(\XX\times\AA_1)$-marginals of the occupation measures induced by the policies
$(\pi_1,\pi_2)\in\bcal{Y}_1\times\bcal{Y}_2$ ---as $\pi_1$ varies and $\pi_2$ remains fixed---
which satisfy the constraint of player 1, that is,
\begin{eqnarray}
\mathcal{A}_{\eta,\rho_1,1}(\pi_{2}) &=& \big\{\mu_{\eta,\pi_1,\pi_2}^{\XX\times\AA_1}: \ \hbox{for $\pi_1\in\bcal{Y}_1$ such that}\ C_1(\eta,\pi_1,\pi_2)\ge\rho_{1}\big\}\nonumber\\
&=&
  \Big\{\gamma\in\mathcal{L}_{\eta,1}(\pi_{2}) :  \int_{\XX\times\AA_1} c^1_1d\gamma +\int_{\XX\times\AA_2} c_1^2d(\gamma^\XX\otimes\pi_2)\geq \rho_{1} \Big\}\subseteq\mathcal{D}_{\eta,1}.\label{eq-def-A1pi2}
\end{eqnarray}
Similarly, for  $\eta\in\bcal{P}^{e}_{\lambda}(\XX)$, $\rho_{2}\in\RR^{p}$, and $\pi_1\in\bcal{Y}_1$ we define
\begin{eqnarray*}
\mathcal{A}_{\eta,\rho_2,2}(\pi_{1}) &= & \big\{\mu_{\eta,\pi_1,\pi_2}^{\XX\times\AA_2}:  \ \hbox{for $\pi_2\in\bcal{Y}_2$ such that}\ C_2(\eta,\pi_1,\pi_2)\ge\rho_2\big\}\\
&=&\Big\{\gamma\in\mathcal{L}_{\eta,2}(\pi_{1}) :\int_{\XX\times\AA_1} c_2^1d(\gamma^\XX\otimes \pi_1)+\int_{\XX\times\AA_2} c_2^2d\gamma \ge\rho_2 \Big\}\subseteq\mathcal{D}_{\eta,2},
\end{eqnarray*}
which is the set of $(\XX\times\AA_2)$-marginals of the occupation measures of the policies $(\pi_1,\pi_2)\in\bcal{Y}_1\times\bcal{Y}_2$ ---as $\pi_1$ is fixed and $\pi_2$ varies---
which satisfy the constraint of player 2. 
We will need the following result in which we will consider an initial distribution $\zeta\in\bcal{P}_\lambda(\XX)$, but not necessarily in $\bcal{P}_\lambda^e(\XX)$.

\begin{proposition}
\label{Contrainte-lower-continuity}
Consider an initial distribution $\zeta\in\bcal{P}_\lambda(\XX)$ and constraint constants $(\rho_{1},\rho_{2})\in\RR^{p}\times\RR^{p}$ that satisfy the Slater condition, and let  $\{\eta_{n}\}_{n\in\NN}$ be an arbitrary sequence in $\bcal{P}^{e}_{\lambda}(\XX)$ converging to  $\zeta\in\bcal{P}_{\lambda}(\XX)$ in total variation.
Then the following assertions hold.
\begin{itemize}
\item[(i).] Fix an arbitrary pair $(\pi_1,\pi_2)\in\bcal{Y}_1\times \bcal{Y}_2$ such that $C_1(\zeta,\pi_1,\pi_2)\ge\rho_1$.
For any sequence 
 $\{\pi_{2,n}\}_{n\in\NN}\subseteq\bcal{Y}_2$ converging to $\pi_2$, there exists a sequence $\{\gamma_n\}_{n\in\NN}$ in $\bcal{P}(\XX\times\AA_1)$ such that $\gamma_n\rightarrow\mu^{\XX\times\AA_1}_{\zeta,\pi_1,\pi_2}$ and such that, for some  $N\in\NN$, we have
  $\gamma_n\in\mathcal{A}_{\eta_n,\rho_1,1}(\pi_{2,n})$ for every $n\ge N$.
\item[(ii).] Fix an arbitrary pair $(\pi_1,\pi_2)\in\bcal{Y}_1\times \bcal{Y}_2$ such that $C_2(\zeta,\pi_1,\pi_2)\ge\rho_2$.
For any sequence 
 $\{\pi_{1,n}\}_{n\in\NN}\subseteq\bcal{Y}_1$ converging to $\pi_1$, there exists a sequence $\{\gamma_n\}_{n\in\NN}$ in $\bcal{P}(\XX\times\AA_2)$ such that $\gamma_n\rightarrow\mu^{\XX\times\AA_2}_{\zeta,\pi_1,\pi_2}$ and such that, for some $N\in\NN$, we have
  $\gamma_n\in\mathcal{A}_{\eta_n,\rho_2,2}(\pi_{1,n})$ for every $n\ge N$.
\end{itemize}
\end{proposition}
\textbf{Proof:} We prove $(i)$. 
From Proposition \ref{Convergence-mu} we have  
$\lim_{n\rightarrow\infty} C_1(\eta_{n},\pi_1,\pi_{2,n})= C_1(\zeta,\pi_1,\pi_2)\ge\rho_{1}$.
Therefore, there exist a sequence $\{\epsilon_{n}\}_{n\in \NN}$ taking values in the interval $(0,1)$ with $\epsilon_n\rightarrow0$ and some index $n_0$ for which
\begin{equation}\label{eq-eq-tool-8}
C_1(\eta_{n},\pi_1,\pi_{2,n})\ge\rho_1-\epsilon_n\quad\hbox{for all $n\ge n_0$}.
\end{equation}
On the other hand, since $\zeta$ and $(\rho_{1},\rho_{2})$ satisfy the Slater condition,
 we can find $\bar\pi_1\in\bcal{Y}_1$ and $\delta>0$ such that 
$C_1(\zeta,\bar\pi_1,\pi_2)>\rho_{1}+\delta$.
So, again from Proposition \ref{Convergence-mu}, 
there is some $n_1\ge n_0$ such that
\begin{equation}\label{eq-eq-tool-9}
C_1(\eta_{n},\bar\pi_1,\pi_{2,n})>\rho_{1}+\delta\quad\hbox{for all $n\ge n_1$}.
\end{equation}
Observe  that for any $n\in\NN$  both $\mu^{\XX\times\AA_1}_{\eta_{n},\pi_1,\pi_{2,n}}$ and $\mu^{\XX\times\AA_1}_{\eta_{n},\bar\pi_1,\pi_{2,n}}$ are in
$\mathcal{L}_{\eta_{n},1}(\pi_{2,n})$, which is a convex set (see Proposition \ref{prop-L1-convex-compact}).
Hence, there is some $\gamma_n\in\mathcal{L}_{\eta_n,1}(\pi_{2,n})$ and  $\pi_{1,n}\in\bcal{Y}_1$ such that 
\begin{align}
\gamma_n=\mu^{\XX\times\AA_1}_{\eta_{n},\pi_{1,n} ,\pi_{2,n}} = (1-\sqrt{\epsilon_{n}}) \mu^{\XX\times\AA_1}_{\eta_{n},\pi_{1} ,\pi_{2,n}}
+ \sqrt{\epsilon_{n}} \mu^{\XX\times\AA_1}_{\eta_{n},\bar{\pi}_{1} ,\pi_{2,n}}
\label{eq-Occup-measure-def-pi1} 
\end{align}
with, as a consequence,
$\mu^{\XX}_{\eta_{n},\pi_{1,n} ,\pi_{2,n}} = (1-\sqrt{\epsilon_{n}}) \mu^{\XX}_{\eta_{n},\pi_{1} ,\pi_{2,n}} + \sqrt{\epsilon_{n}} \mu^{\XX}_{\eta_{n},\bar{\pi}_{1} ,\pi_{2,n}}$,
and so
\begin{eqnarray}
\mu^{\XX\times\AA_2}_{\eta_{n},\pi_{1,n} ,\pi_{2,n}} &=& \mu_{\eta_{n},\pi_{1,n},\pi_{2,n}}^\XX\otimes \pi_{2,n}
=(1-\sqrt{\epsilon_{n}}) \mu^{\XX\times\AA_2}_{\eta_{n},\pi_{1} ,\pi_{2,n}} + \sqrt{\epsilon_{n}} \mu^{\XX\times\AA_2}_{\eta_{n},\bar{\pi}_{1} ,\pi_{2,n}}.
\label{eq-eq-tool-19}
\end{eqnarray}
It is clear from \eqref{eq-Occup-measure-def-pi1} and Proposition \ref{Convergence-mu} that 
$\lim_{n\rightarrow\infty} \gamma_n=\lim_{n\rightarrow\infty} \mu^{\XX\times\AA_1}_{\eta_n,\pi_1,\pi_{2,n}}= \mu_{\zeta,\pi_1,\pi_2}^{\XX\times\AA_1}.$
We deduce from \eqref{eq-Occup-measure-def-pi1}--\eqref{eq-eq-tool-19} that for all $n\ge n_1$
\begin{eqnarray}
C_1(\eta_{n},\pi_{1,n},\pi_{2,n}) &=& \int_{\XX\times\AA_1} c^1_1 d\mu_{\eta_{n},\pi_{1,n},\pi_{2,n}}^{\XX\times\AA_1}+
\int_{\XX\times\AA_1} c^2_1 d\mu_{\eta_{n},\pi_{1,n},\pi_{2,n}}^{\XX\times\AA_2}\nonumber\\
&=& (1-\sqrt{\epsilon_n})C_1(\eta_{n},\pi_1,\pi_{2,n})+\sqrt{\epsilon_n}C_1(\eta_{n},\bar\pi_1,\pi_{2,n})\nonumber\\  \noalign{\smallskip}
&\ge&
\rho_{1} + \sqrt{\epsilon_{n}} \big[\delta -(1-\sqrt{\epsilon_{n}}) \sqrt{\epsilon_{n}} \big], \label{eq-eq-tool-18}
\label{eq-Occup-measure-def-pi2}
\end{eqnarray}
where the last inequality is derived from \eqref{eq-eq-tool-8} and \eqref{eq-eq-tool-9}. Consequently,  
there exists some $N\ge n_1$ such that  $C_1(\eta_{n},\pi_{1,n},\pi_{2,n})\ge\rho_{1}$ for $n\ge N$. Since by definition $\gamma_n=\mu^{\XX\times\AA_1}_{\eta_{n},\pi_{1,n} ,\pi_{2,n}}$,
this establishes precisely that $\gamma_n\in\mathcal{A}_{\eta_n,\rho_1,1}(\pi_{2,n})$ for all $n\ge N$. This completes the proof of (i).
The rationale for using  the coefficient $\sqrt{\epsilon_n}$ in \eqref{eq-Occup-measure-def-pi1}  is to mix the measures at a rate slower than the bound in
\eqref{eq-eq-tool-8} in order to satisfy the constraint as in \eqref{eq-eq-tool-18}.
\hfill $\Box$

\begin{proposition}
\label{Continuity-A-1} Let  $\eta\in \bcal{P}^{e}_{\lambda}(\XX)$ and $(\rho_{1},\rho_{2})\in \RR^{p}\times\RR^{p}$
satisfy the Slater condition. Then, 
the correspondences $\mathcal{A}_{\eta,\rho_1,1}:\bcal{Y}_2\tto\mathcal{D}_{\eta,1}$ defined by $\pi_{2}\mapsto \mathcal{A}_{\eta,\rho_1,1}(\pi_{2})$ and 
$\mathcal{A}_{\eta,\rho_2,2}:\bcal{Y}_1\tto\mathcal{D}_{\eta,2}$ given  by $\pi_{1}\mapsto \mathcal{A}_{\eta,\rho_2,2}(\pi_{1})$ are both continuous with nonempty convex and compact values.
\end{proposition}
\textbf{Proof.} We make the proof for the correspondence $\mathcal{A}_{\eta,\rho_1,1}$.
Given $\pi_{2}\in \bcal{Y}_{2}$, we have that  $\mathcal{A}_{\eta,\rho_1,1}(\pi_{2})$ is nonempty by Assumption \ref{Assump-Slater}.
Since $\mathcal{L}_{\eta,1}(\pi_2)$ is convex (recall  Proposition \ref{prop-L1-convex-compact}), the fact that $\mathcal{A}_{\eta,\rho_1,1}(\pi_2)$ is given by a linear constraint in
$\gamma$ ---see \eqref{eq-def-A1pi2}--- yields that it is indeed a convex set.

We prove compactness of $\mathcal{A}_{\eta,\rho_1,1}(\pi_{2})$ along with upper semicontinuity of $\pi_2\mapsto\mathcal{A}_{\eta,\rho_1,1}(\pi_2)$. To this end, we use the Closed Graph Theorem in \cite[Theorem 17.11]{aliprantis06}. Indeed, since the range of the correspondence is the compact metric space $\mathcal{D}_{\eta,1}$, the correspondence has closed graph if and only if it is upper semicontinuous and closed-valued (hence, in our case, compact-valued). So, it suffices to show that $\pi_2\mapsto\mathcal{A}_{\eta,\rho_1,1}(\pi_2)$ has closed graph.
Suppose that we have a convergent sequence $(\pi_{2,n},\gamma_n)$ in the graph of $\mathcal{A}_{\eta,\rho_1,1}$. This means that  $\pi_{2,n}\rightarrow \pi_{2}$ in $\bcal{Y}_{2}$
and that $\gamma_{n}\in\mathcal{A}_{\eta,\rho_1,1}(\pi_{2,n})$ are such that $\gamma_{n}\rightarrow \gamma\in\mathcal{D}_{\eta,1}$. Our goal is to prove that
$\gamma\in\mathcal{A}_{\eta,\rho_1,1}(\pi_2)$. We note that  for each $n\in\NN$ there exists $\pi_{1,n}\in \bcal{Y}_{1}$ such that
\begin{align}
\gamma_{n} = \mu_{\eta,\pi_{1,n},\pi_{2,n}}^{\XX\times\AA_{1}} \quad\hbox{and}\quad C_1(\eta,\pi_{1,n},\pi_{2,n})\ge\rho_1
\label{eq-tool-again}
\end{align}
for all $n\in\NN$.
There exists a subsequence $\{\pi_{1,n'}\}$ of $\{\pi_{1,n}\}_{n\in\NN}$ that converges to some $\pi_{1}\in\bcal{Y}_{1}$.
Using Proposition \ref{Convergence-mu} in \eqref{eq-tool-again}, we have 
$\gamma_{n'}\rightarrow \mu_{\eta,\pi_{1},\pi_{2}}^{\XX\times\AA_{1}}=\gamma$,
which establishes $\gamma\in \mathcal{L}_{\eta,1}(\pi_{2})$,
and 
$C_1(\eta,\pi_{1},\pi_{2})\ge\rho_1$, from which $\gamma\in\mathcal{A}_{\eta,\rho_1,1}(\pi_2)$ follows.
Lower semicontinuity of $\mathcal{A}_{\eta,\rho_1,1}$ is a direct consequence of Proposition \ref{Contrainte-lower-continuity}(i) applied to the constant sequence $\eta_n\equiv\eta$ 
and  the sequential characterization of lower semicontinuity of correspondences given in  \cite[Theorem 17.21]{aliprantis06}. 
\hfill $\Box$
\\[10pt]\indent
Based on Remark \ref{remark-cost}, given $\pi_1\in\bcal{Y}_1$ or, equivalently, given
$\gamma\in\mathcal{L}_{\eta,1}(\pi_2)$ we have
$$R_1(\eta,\pi_1,\pi_2)=\int_{\XX\times\AA_1} r^1_1d\gamma+\int_{\XX\times\AA_2} r^2_1 d(\gamma^\XX\otimes\pi_2).$$
For a fixed policy of player 2, therefore, the goal of player 1 is to maximize  the above integral over all $\gamma\in\mathcal{A}_{\eta,\rho_1,1}(\pi_2)$, that is, over all measures $\gamma\in\mathcal{L}_{\eta,1}(\pi_2)$ which satisfy his constraint. 
This leads  to the definition of the correspondence $\mathcal{H}_{\eta,\rho_1,1}:\bcal{Y}_2\tto \mathcal{D}_{\eta,1}$ given by
\begin{eqnarray*}
\mathcal{H}_{\eta,\rho_{1},1}(\pi_{2}) 
=\argmax_{\gamma\in\mathcal{A}_{\eta,\rho_1,1}(\pi_{2})} \Bigg\{\int_{\XX\times\AA_1} r^1_1d\gamma+ \int_{\XX\times\AA_2} r^2_1d(\gamma^{\XX}\otimes\pi_{2})\Bigg\}
\end{eqnarray*}
for $\pi_2\in\bcal{Y}_2$, $\rho_{1}\in\RR^{p}$, and $\eta\in \bcal{P}^{e}_{\lambda}(\XX)$.
Similarly, we define  $\mathcal{H}_{\eta,\rho_2,2}:\bcal{Y}_1\tto \mathcal{D}_2$ 
as
\begin{eqnarray*}
\mathcal{H}_{\eta,\rho_{2},2}(\pi_{1}) 
&=&\argmax_{\gamma\in\mathcal{A}_{\eta,\rho_2,2}(\pi_{1})} \Bigg\{\int_{\XX\times\AA_1} r^1_2 d(\gamma^{\XX}\otimes\pi_{1}) + \int_{\XX\times\AA_2} r^2_2d\gamma\Bigg\}
\end{eqnarray*}
for $\pi_1\in\bcal{Y}_1$, $\rho_{2}\in\RR^{p}$, and $\eta\in \bcal{P}^{e}_{\lambda}(\XX)$

\begin{proposition}
\label{Correspondence-H1}
Let $\eta\in \bcal{P}^{e}_{\lambda}(\XX)$ and $(\rho_{1},\rho_{2})\in \RR^{p}\times\RR^{p}$
satisfy the Slater condition.
The correspondences $\mathcal{H}_{\eta,\rho_{1},1}$ and $\mathcal{H}_{\eta,\rho_{2},2}$ are upper semicontinuous with nonempty compact and convex values.
\end{proposition}
\textbf{Proof.} We prove the result only for $\mathcal{H}_{\eta,\rho_{1}1}$.
We have that 
$$
(\pi_2,\gamma)\mapsto f(\pi_2,\gamma)= \int_{\XX\times\AA_1} r^1_1d\gamma+ \int_{\XX\times\AA_2} r^2_1d(\gamma^{\XX}\otimes\pi_{2})
$$
is continuous on the graph of $\mathcal{A}_{\eta,\rho_1,1}$. Indeed, suppose that $\{\pi_{2,n}\}\subseteq\bcal{Y}_2$ converges to some $\pi_2\in\bcal{Y}_2$ and that
$\gamma_n\in\mathcal{A}_{\eta,\rho_1,1}(\pi_{2,n})$ is such that $\gamma_n\rightarrow\gamma$ with, necessarily (the correspondence $\mathcal{A}_{\eta,\rho_1,1}$ being closed) $\gamma\in\mathcal{A}_{\eta,\rho_1,1}(\pi_2)$. First of all note that $r_1^1\in\car_b(\XX\times\AA_1,\RR)$ and so by definition of the $ws$-convergence  we have 
$ \int_{\XX\times\AA_1} r^1_1d\gamma_n\rightarrow \int_{\XX\times\AA_1}r^1_1d\gamma$.
On the other hand, since the sequence $\int r^2_1 d(\gamma_n^\XX\otimes\pi_{2,n})$ is bounded, in order to prove the convergence
\begin{equation}\label{eq-tool-20}
\int_{\XX\times\AA_2} r^2_1 d(\gamma_n^\XX\otimes\pi_{2,n})\rightarrow
\int_{\XX\times\AA_2} r^2_1 d(\gamma^\XX\otimes\pi_{2})
\end{equation}
it suffices to show that the above limit holds through any convergent subsequence $\int r^2_1 d(\gamma_{n'}^\XX\otimes\pi_{2,n'})$. For $\gamma_{n'}$ there exists some $\pi_{1,n'}\in\bcal{Y}_1$ such that $\gamma_{n'}=\mu^{\XX\times\AA_1}_{\eta,\pi_{1,n'},\pi_{2,n'}}$.
Without loss of generality we can assume that, for some $\pi_1\in\bcal{Y}_1$, we have $\pi_{1,n'}\rightarrow\pi_1$. In particular, using Proposition \ref{Convergence-mu}, this implies that $\gamma_{n'}\rightarrow \mu^{\XX\times\AA_1}_{\eta,\pi_1,\pi_2}=\gamma$. But we also have
$\gamma_{n'}^\XX\otimes\pi_{2,n'}=\mu^{\XX\times\AA_2}_{\eta,\pi_{1,n'},\pi_{2,n'}}\rightarrow\mu^{\XX\times\AA_2}_{\eta,\pi_{1},\pi_{2}}=\gamma^\XX\otimes\pi_2$
by using again Proposition~\ref{Convergence-mu}.
Thus, since $r^2_1\in\car_b(\XX\times\AA_2,\RR)$, the convergence \eqref{eq-tool-20} through~$n'$ follows. This completes the proof of the continuity of the function $f$ on the graph of $\mathcal{A}_{\eta,\rho_1,1}$.

Now we are in position to apply  Berge's Maximum Theorem  \cite[Theorem 17.31]{aliprantis06} to the function~$f$ on the graph of $\mathcal{A}_{\eta,\rho_1,1}$. Indeed, we have already shown that $f$ is continuous and, besides, we have that $\mathcal{A}_{\eta,\rho_1,1}$ is a continuous correspondence with nonempty compact values (Proposition \ref{Continuity-A-1}). By the maximum theorem we conclude that the $\argmax$ correspondence  (that is, $\mathcal{H}_{\eta,\rho_1,1}$)   is upper semicontinuous with nonempty compact values.

Finally, from Proposition \ref{Continuity-A-1} we also have that $\mathcal{A}_{\eta,\rho_1,1}(\pi_{2})$ is convex and since the mapping 
$\gamma\mapsto f(\pi_2,\gamma)= \int_{\XX\times\AA_1} r^1_1d\gamma+ \int_{\XX\times\AA_2} r^2_1d(\gamma^{\XX}\otimes\pi_{2})$
is linear  in $\gamma$ for fixed $\pi_2$, it follows that its set of maxima is a convex set.  This proves that $\mathcal{H}_{\eta,\rho_1,1}$ is also convex-valued.
\hfill $\Box$
\\[10pt]\indent
Hence, given an initial distribution $\eta\in \bcal{P}^{e}_{\lambda}(\XX)$ and constraint constants $\rho=(\rho_{1},\rho_{2})\in\RR^{p}\times\RR^{p}$ satisfying the Slater condition, let us define the correspondence
$\mathcal{H}_{\eta,\rho}:\bcal{Y}_{1}\times\bcal{Y}_{2}\tto \mathcal{D}_{\eta,1}\times\mathcal{D}_{\eta,2}$ as
$$
\mathcal{H}_{\eta,\rho}(\pi_{1},\pi_{2})=\mathcal{H}_{\eta,\rho_{1},1}(\pi_{2})\times \mathcal{H}_{\eta,\rho_{2},2}(\pi_{1})
\quad\hbox{for $(\pi_1,\pi_2)\in\bcal{Y}_1\times\bcal{Y}_2$}.
$$
As a consequence of the previous results and using \cite[Theorem 17.28]{aliprantis06}, we conclude that the correspondence $\mathcal{H}_{\eta,\rho}$
is upper semicontinuous with nonempty compact and convex values. 

\paragraph{Nash equilibrium for an initial distribution $\eta\in\bcal{P}_\lambda^e(\XX)$.}

\begin{theorem}\label{th-main-bis}
Let $\eta\in \bcal{P}^{e}_{\lambda}(\XX)$ and $\rho=(\rho_{1},\rho_{2})\in \RR^{p}\times\RR^{p}$ satisfy the Slater condition.  Under Assumption \ref{Assump} the following results hold.
\begin{itemize} 
\item[(i)]  The  correspondence 
$\mathcal{H}_{\eta,\rho}\circ\mathcal{J}_{\eta}:\mathcal{D}_{\eta,1}\times\mathcal{D}_{\eta,2}\tto \mathcal{D}_{\eta,1}\times\mathcal{D}_{\eta,2}$
has a fixed point $(\gamma^*_1,\gamma^*_2)\in\mathcal{D}_{\eta,1}\times\mathcal{D}_{\eta,2}$.
\item[(ii)]
The pair of stationary Markov policies  $(\pi^*_1,\pi^*_2)=\mathcal{J}_{\eta}(\gamma^*_1,\gamma^*_2)\in\bcal{Y}_1\times\bcal{Y}_2$ is a constrained equilibrium in the class of history-dependent policies $\mathbf{\Pi}_1\times\mathbf{\Pi}_2$  for the game  model $\mathcal{G}(\eta,\rho)$.\end{itemize}
\end{theorem}
\textbf{Proof:} (i). We can use \cite[Theorem 17.23]{aliprantis06} to prove that  $\mathcal{H}_{\eta,\rho}\circ\mathcal{J}_{\eta}$ is an upper semicontinuous correspondence (it is the composition of the continuous function $\mathcal{J}_{\eta}$ and the upper semicontinuous correspondence $\mathcal{H}_{\eta,\rho}$)  with nonempty compact and convex values.  By the Closed Graph Theorem \cite[Theorem 17.11]{aliprantis06}, the correspondence $\mathcal{H}_{\eta,\rho}\circ\mathcal{J}_{\eta}$ is closed.
On the other hand,  $\mathcal{D}_{\eta,1}\times\mathcal{D}_{\eta,2}$ is a nonempty compact convex subset  of the locally convex Hausdorff space
$\bcal{M}(\XX\times\AA_{1})\times\bcal{M}(\XX\times\AA_{2})$; recall Proposition \ref{prop-locally-convex}.
The Kakutani-Fan-Glicksberg fixed point theorem \cite[Corollary 17.55]{aliprantis06} yields the existence of a fixed point for $\mathcal{H}_{\eta,\rho}\circ\mathcal{J}_{\eta}$.\\[5pt]
\noindent
(ii). Denote by $(\gamma^*_1,\gamma^*_2)\in\mathcal{D}_{\eta,1}\times\mathcal{D}_{\eta,2}$  a fixed point of 
$\mathcal{H}_{\eta,\rho}\circ\mathcal{J}_{\eta}$ and let $(\pi^*_1,\pi^*_2)\in\bcal{Y}_1\times\bcal{Y}_2$ be the stationary Markov policies given by
$(\pi^*_1,\pi^*_{2})=\mathcal{J}_{\eta}(\gamma^*_1,\gamma^*_2)$,
so that 
$(\gamma^*_1,\gamma^*_2)\in\mathcal{H}_{\eta,\rho}(\pi^*_1,\pi^*_{2})$.
Observe that $\gamma^*_1\in\mathcal{L}_{\eta,1}(\pi^*_2)$ and so $\gamma^*_1=\mu^{\XX\times\AA_1}_{\eta,\pi_1,\pi^*_2}$ for some
$\pi_1\in\bcal{Y}_1$. But then we have $\pi_1=\mathcal{J}_{\eta,1}(\gamma^*_1)=\pi^*_1$. Hence, using the symmetric argument for $\gamma^*_2$ we can conclude that 
\begin{equation}\label{eq-tool-21}
\gamma^*_1=\mu^{\XX\times\AA_1}_{\eta,\pi^*_1,\pi^*_2}
\quad\hbox{and}\quad
\gamma^*_2=\mu^{\XX\times\AA_1}_{\eta,\pi^*_1,\pi^*_2}.
\end{equation}
It follows that
$ C_1(\eta,\pi^*_1,\pi^*_2)=\int_{\XX\times\AA_1} c^1_1d\gamma^*_1+\int_{\XX\times\AA_2} c^2_1 d(\gamma_1^{*\XX}\otimes\pi^*_2)
 \ge\rho_{1}$
since $\gamma^*_1\in\mathcal{A}_{\eta,\rho_1,1}(\pi^*_2)$.
Similarly, we can show that $C_2(\eta,\pi^*_1,\pi^*_2)\ge\rho_2$ and so,
 the stationary Markov policies $(\pi^*_1,\pi^*_2)$
satisfy the constraints of both players (see item (i) in Definition \ref{def-nash}).

Suppose now that player $1$ varies his policy from $\pi^*_1\in\bcal{Y}_1$ to some history-dependent policy $\pi_1\in\mathbf{\Pi}_1$ which satisfies his own constraint, that is, $C_1(\eta,\pi_1,\pi^*_2)\ge\rho_{1}$. We want to check that $R_1(\eta,\pi_1,\pi^*_2)\le R_1(\eta,\pi^*_1,\pi^*_2)$. 
By Corollary \ref{cor-Markov-sufficient} we can suppose without loss of generality that $\pi_1\in\bcal{Y}_1$. Hence, for this Markov policy $\pi_1\in\bcal{Y}_1$ we have
\begin{equation}\label{eq-tool-23}
C_1(\eta,\pi_1,\pi_2^*)=\int_{\XX\times\AA_1} c^1_1d\mu^{\XX\times\AA_1}_{\eta,\pi_1,\pi^*_2}
+\int_{\XX\times\AA_2} c^2_1d\mu^{\XX\times\AA_2}_{\eta,\pi_1,\pi^*_2}
\ge\rho_{1}.
\end{equation}
With this in mind, \eqref{eq-tool-23}
 implies that $\mu^{\XX\times\AA_1}_{\eta,\pi_1,\pi^*_2}\in\mathcal{A}_{\eta,\rho_1,1}(\pi^*_2)$ but by the definition of
 $\gamma^*_1\in\mathcal{H}_{\eta,\rho_{1},1}(\pi^*_2)$ we have
\begin{eqnarray*}
R_1(\eta,\pi_1,\pi_2^*) &=&
\int_{\XX\times\AA_1} r^1_1d\mu^{\XX\times\AA_1}_{\eta,\pi_1,\pi^*_2}+ \int_{\XX\times\AA_2} r^2_1d(\mu^{\XX}_{\eta,\pi_1,\pi^*_2}\otimes\pi^*_{2})\\
&\le&
\int_{\XX\times\AA_1} r^1_1d\gamma^*_1+ \int_{\XX\times\AA_2} r^2_1d(\gamma_1^{*\XX}\otimes\pi^*_{2})
= R_1(\eta,\pi^*_1,\pi^*_2),
\end{eqnarray*}
where we have made use of \eqref{eq-tool-21}.
We proceed similarly for player 2 and we conclude that the stationary Markov policies $(\pi^*_1,\pi^*_2)\in\bcal{Y}_1\times\bcal{Y}_2$ are indeed a constrained Nash equilibrium. 
\hfill$\Box$

\subsection{Existence of a Nash equilibrium for an initial distribution $\nu\in\bcal{P}_\lambda(\XX)$}
\label{section42}

Theorem \ref{th-main-bis}   establishes the existence of a Nash equilibrium for an initial distribution in $\mathcal{P}_\lambda^e(\XX)$. To obtain our main result in this paper, we now drop this condition and replace it with $\nu\in\mathcal{P}_\lambda(\XX)$.
\begin{theorem}\label{th-main}
Consider the game model $\mathcal{G}(\nu,\theta)$ 
and suppose that Assumption \ref{Assump} is satisfied for the initial distribution $\nu\in\bcal{P}_\lambda(\XX)$. 
There exist stationary Markov policies $(\pi^*_1,\pi^*_2)\in\mathbf{M}_1\times\mathbf{M}_2$ which are a constrained equilibrium for the players in the class of history-dependent policies $\mathbf{\Pi}_1\times\mathbf{\Pi}_2$.
\end{theorem}
\textbf{Proof.} Let $\{\eta_{n}\}_{n\in\NN}\subseteq\bcal{P}_\lambda^e(\XX)$ be the sequence defined by $\eta_{n}=\frac{n}{n+1}\nu+\frac{1}{n+1}\lambda$ for $n\in\NN$. 
Clearly,  $\{\eta_{n}\}_{n\in\NN}$ converges to $\nu$ in total variation.
Recalling \eqref{Def-marginal-occup-measure}--\eqref{eq-occupation-marginal-conditional} and since the constraint function $c_i$ is bounded,
there is some constant $\mathfrak{c}>0$ such that 
$|C_i(\eta_{n},\pi_{1},\pi_{2})-C_i(\nu,\pi_{1},\pi_{2})| \leq \mathfrak{c}/(n+1)$ for all $i=1,2$ and $(\pi_{1},\pi_{2})\in\bcal{Y}_{1}\times\bcal{Y}_{2}$.
Hence, the initial distribution $\eta_{n}\in\bcal{P}_\lambda^e(\XX)$ and the constraint constants $\rho_n=(\rho_{1,n},\rho_{2,n})$ given by
$\rho_{i,n}=\theta_{i}-\frac{\mathfrak{c}}{n+1}\mathbf{1}_p$ for $n\in\NN$
satisfy the Slater condition in Definition \ref{def-slater}. We are thus in position to apply
Theorem \ref{th-main-bis} to the game model $\mathcal{G}(\eta_{n},\rho_{n})$ to obtain the existence of a Nash equilibrium. That is, for each $n\in\NN$ there exist
stationary Markov policies $\{(\pi^*_{1,n},\pi^*_{2,n})\}_{n\in\NN}$ in $\bcal{Y}_1\times\bcal{Y}_2$ 
satisfying 
\begin{align}
\label{ineq2-sensitivity}
C_i(\eta_{n},\pi^*_{1,n},\pi^*_{2,n})\ge\rho_{i,n}\quad\hbox{for $i=1,2$.}
\end{align}
and
$$\forall\pi_1\in\mathbf{\Pi}_1,\ C_1(\eta_{n},\pi_{1},\pi^*_{2,n})\ge\rho_{1,n}  \Rightarrow  R_1(\eta_{n},\pi^*_{1,n},\pi^*_{2,n})\ge R_1(\eta_{n},\pi_1,\pi^*_{2,n}),$$
$$\forall\pi_2\in\mathbf{\Pi}_2,\ C_2(\eta_{n},\pi^*_{1,n},\pi_{2})\ge\rho_{2,n}  \Rightarrow  R_2(\eta_{n},\pi^*_{1,n},\pi^*_{2,n})\ge R_2(\eta_{n},\pi^*_{1,n},\pi_2).$$
By compactness of $\bcal{Y}_1$ and $\bcal{Y}_2$, there exists a convergent subsequence of $\{(\pi^*_{1,n},\pi^*_{2,n})\}_{n\in\NN}$. Without loss of generality, we will assume that the whole sequence converges   to some 
$(\pi^*_{1},\pi^*_{2})\in\bcal{Y}_1\times\bcal{Y}_2$. Our goal now is to show that $(\pi^*_{1},\pi^*_{2})$ is a Nash equilibrium for the game model $\mathcal{G}(\nu,\theta)$.

From Proposition \ref{Convergence-mu}, we get by taking the limit in \eqref{ineq2-sensitivity}
$C_i(\nu,\pi^*_{1},\pi^*_{2})\ge\theta_{i}$ for $i=1,2$.
On the other hand, 
let $\pi_1$ be an arbitrary policy in $\mathbf{\Pi}_{1}$ satisfying $C_1(\nu,\pi_{1},\pi^*_{2})\ge\theta_{1}$. Without loss of generality, it can be assumed that
$\pi_{1}\in \bcal{Y}_{1}$ (recall Corollary \ref{cor-Markov-sufficient}). Consider the associated measure $\gamma=\mu_{\nu,\pi_1,\pi^*_2}^{\XX\times\AA_1}$.
We can apply Proposition \ref{Contrainte-lower-continuity} to deduce
the existence of a sequence $\gamma_n\rightarrow\gamma$ with $\gamma_n\in\mathcal{A}_{\eta_n,\theta_{1},1}(\pi^*_{2,n})$ for large enough $n\ge N$.
For such large $n$, let $\pi_{1,n}\in\bcal{Y}_1$ be such that  $\gamma_n=\mu^{\XX\times\AA_1}_{\eta_n,\pi_{1,n},\pi^*_{2,n}}$.
We have that
$C_1(\eta_{n},\pi_{1,n},\pi^{*}_{2,n})\ge\theta_1\ge\rho_{1,n}$ for any $n\geq N$. But $(\pi^*_{1,n},\pi^*_{2,n})$ being a Nash equilibrium for the game model $\mathcal{G}(\eta_n,\rho_{n})$,    for such $n\geq N$
\begin{eqnarray*}
R_1(\eta_{n},\pi^*_{1,n},\pi^*_{2,n}) &\ge&  R_1(\eta_{n},\pi_{1,n},\pi^*_{2,n}) = \int_{\XX\times\AA_1} r_1^1 d\gamma_n +\int_{\XX\times\AA_2} r_1^2 d(\gamma_n^\XX\otimes\pi^*_{2,n}).
\end{eqnarray*}
Observe that since $\gamma_n\rightarrow\gamma$, we have 
$$ \lim_{n\rightarrow\infty}  \int_{\XX\times\AA_1} r_1^1 d\gamma_n=
\int_{\XX\times\AA_1} r^1_1 d\gamma=
\int_{\XX\times\AA_1} r^1_1 d \mu^{\XX\times\AA_1}_{\nu,\pi_1,\pi^*_2}.$$ 
We want to show that
$ \int_{\XX\times\AA_2} r_1^2 d(\gamma_n^\XX\otimes\pi^*_{2,n})\rightarrow \int_{\XX\times\AA_2} r^2_1 d\mu^{\XX\times\AA_2}_{\nu,\pi_1,\pi^*_2}$. Since $r^2_1$ is bounded, it suffices to show that the limit holds through any convergent subsequence  $\{n'\}$ of $\int_{\XX\times\AA_2} r_1^2 d(\gamma_n^\XX\otimes\pi^*_{2,n})$. Without loss of generality, we can assume that $\pi_{1,n'}\rightarrow\bar\pi_1$ for some $\bar\pi_1\in\bcal{Y}_1$.
By Proposition \ref{Convergence-mu} we have
$$\gamma_{n'}^\XX\otimes\pi^*_{2,n'}=\mu^{\XX\times\AA_2}_{\eta_{n'},\pi_{1,n'},\pi^*_{2,n'}}\rightarrow\mu^{\XX\times\AA_2}_{\nu,\bar\pi_{1},\pi^*_{2}}=\mu^\XX_{\nu,\bar\pi_1,\pi^*_2}\otimes\pi^*_2$$ and, in particular, $\gamma_{n'}^\XX\rightarrow \mu^{\XX}_{\nu,\bar\pi_1,\pi^*_2}$ in the $s$-topology of $\bcal{P}(\XX)$.
On the other hand, we have $\gamma_{n'}^\XX\rightarrow \gamma^\XX=\mu^\XX_{\nu,\pi_1,\pi^*_2}$ is the $s$-topology of $\bcal{P}(\XX)$.
This implies that $\mu^{\XX}_{\nu,\pi_1,\pi^*_2}=\mu^\XX_{\nu,\bar\pi_1,\pi^*_2}$. 
We conclude that 
$$\gamma_{n'}^\XX\otimes\pi^*_{2,n'}\rightarrow \mu^\XX_{\nu,\pi_1,\pi^*_2}\otimes\pi^*_2=\mu^{\XX\times\AA_2}_{\nu,\pi_1,\pi^*_2}$$
and the desired convergence follows.
Using again Proposition \ref{Convergence-mu} we can take the limit in
$$
R_1(\eta_{n},\pi^*_{1,n},\pi^*_{2,n}) \ge   \int_{\XX\times\AA_1} r_1^1 d\gamma_n +\int_{\XX\times\AA_2} r_1^2 d(\gamma_n^\XX\otimes\pi^*_{2,n}).
$$
and obtain
$R_1(\nu,\pi^*_1,\pi^*_2)\ge R_1(\nu,\pi_1,\pi^*_2)$.
Summarizing, if $\pi_1\in\bcal{Y}_1$ satisfies $C_1(\nu,\pi_1,\pi^*_2)\ge\theta_1$ then, necessarily, we have $R_1(\nu,\pi^*_1,\pi^*_2)\ge R_1(\nu,\pi_1,\pi^*_2)$.

Proceeding symmetrically for player 2, it follows that $(\pi^*_1,\pi^*_2)$ is a Nash equilibrium.
\hfill$\Box$

\end{document}